\documentclass[smallextended,referee,envcountsect,]{svjour3}
\smartqed
\usepackage[square,sort,comma,numbers]{natbib}
\usepackage{amsmath}
\usepackage{amssymb}
\usepackage{graphicx}
\usepackage{xcolor}
\journalname{arXiv}

\begin{document}

\title{A Complete Inverse Optimality Study for a Tank-Liquid System}

\author{Iasson Karafyllis \and Filippos Vokos \and Miroslav Krstic }

\institute{Iasson Karafyllis,  Corresponding author  \\
iasonkar@central.ntua.gr \\
Filippos Vokos \\
fivojean@mail.ntua.gr \at
             Dept. of Mathematics, National Technical University of Athens \\
              Zografou Campus, 15780, Athens, Greece\\
           \and        
              Miroslav Krstic  \at
              krstic@ucsd.edu \\
              Dept. of Mechanical and Aerospace Eng., University of California \\
              San Diego, La Jolla, CA 92093-0411, U.S.A. \\
}


\maketitle

\begin{abstract}
This paper presents a complete inverse optimality study for a linearized tank-liquid system where the liquid is described by the viscous Saint-Venant model with surface tension and possible wall friction. We define an appropriate weak solution notion for which we establish existence/uniqueness results with inputs that do not necessarily satisfy any compatibility condition as well as stabilization results with feedback laws that are constructed with the help of a Control Lyapunov Functional. We show that the proposed family of stabilizing feedback laws is optimal for a certain meaningful quadratic cost functional. Finally, we show that the optimal feedback law guarantees additional stronger stability estimates which are similar to those obtained in the case of classical solutions. 
\end{abstract}
\keywords{Saint-Venant \and inverse optimality \and Control Lyapunov Functional \and PDEs \and feedback stabilization}
\subclass{93C20 \and  93B52 \and 49N45 \and 35G16}


\section{Introduction}

The concept of inverse optimality, i.e., the study of whether a feedback law is optimal for certain meaningful cost functionals, is crucial in control theory. 
Inverse optimality is a property of a feedback law where, instead of seeking a feedback law which is optimal for a given cost criterion, which is an essentially unsolvable problem for nonlinear systems with more than a few states due to the curse of dimensionality, a family of feedback laws is explicitly constructed instead, which ensures that some cost functional, which is meaningful though not given a priori, is minimized in closed loop. In the case of linear control systems, a meaningful cost functional is the sum of two quadratic, coercive functionals: one for the state and one for the control input.
This process is (almost) the opposite process of optimal control, where we first select a specific cost functional and then try to find the corresponding optimal control- in the form of a feedback law or not (see for instance \cite{B6,26,34}). 

Inverse optimal control was invented by Rudolf Kalman in his celebrated 
paper which poses the question whether every stabilizing feedback law for a 
linear system is optimal with respect to a meaningful cost functional \cite{11}. Kalman answers 
this question in the negative but characterizes all feedback laws 
that are inverse optimal: those whose gain matrix is the product of the transpose of 
the input matrix and the Lyapunov matrix of the closed-loop stabilized system. 

The extensions from Kalman’s work for linear systems to inverse optimal control for 
nonlinear systems affine in control begins with Moylan and Anderson \cite{28}. Two decades later, 
a complete methodology for the design of robust inverse optimal nonlinear controllers for 
control-affine systems was developed in \cite{9, 23, 25, 31}, as well as the extensions in \cite{6,7} to stochastic systems. More recently, inverse optimality has been successfully extended
to delay systems in \cite{3} and has been adapted to the problem of inverse optimal safety filter design
in \cite{24}.


Despite the importance of inverse optimality in control theory, there are very few works studying global inverse optimality for control systems described by Partial Differential Equations (PDEs) with boundary control. There are reasons that can explain the scarcity of global inverse optimality results. The first reason is the fact that inverse optimality requires a Lyapunov-based feedback design procedure (of LgV type; see Chapter 4 in \cite{9} and Chapter 8 in \cite{33}) which is not frequent in the study of control systems described by PDEs. Indeed, the Control Lyapunov Function (CLF) methodology (see \cite{12, 33} and references therein) is one of the most important methodologies for the construction of stabilizing feedback laws for systems described by Ordinary Differential Equations (ODEs). When studying systems described by PDEs, the Control Lyapunov Function becomes a Control Lyapunov Functional (CLF) and the use of CLFs for the solution of global feedback stabilization problems for systems with PDEs has been presented in detail in \cite{5}. CLFs for PDEs have been used for instance in \cite{1, 10, 13, 15, 22, 27, 29, 30}, although in most of these works the constructed feedback law is not of the type that can be shown to be optimal for some meaningful cost functional (it is not of LgV type). However, it was shown in \cite{32} that certain backstepping procedures can give optimal boundary feedback laws for parabolic PDEs with Neumann actuation. 
In \cite{B2} the inverse optimal adaptive control problem is studied for a parabolic PDE with an unknown reaction coefficient while in \citep{B3} it is proved that the boundary control constructed in \cite{1} for the Korteweg-de Vries-Burgers equation is inverse optimal. Other interesting studies are \cite{B5}, in which the inverse optimal control problem of a class of linear hyperbolic PDEs is considered, and \cite{B4} where an inverse optimal gain assignment problem is solved for evolution systems in Hilbert spaces. A study concerning the construction of an optimal Lyapunov-based boundary feedback law which maximizes a meaningful cost functional is \cite{C1}. In \cite{C1} the feedback controller achieves mixing in a 3D pipe flow described by the Navier-Stokes equations.

Another reason that can explain the scarcity of global inverse optimality results for control systems described by PDEs with boundary control is the difficulty of the technical issues that arise in such problems. Usually, feedback stabilization problems are solved for classical solutions of PDEs. Classical solutions require the so-called compatibility conditions (see \cite{2, 8}) which actually -in the case of linear PDEs- mean stabilizability on a subspace of the state space. This is not really a problem in stabilization problems: the compatibility conditions can be satisfied after an arbitrarily short transient period during which an input that establishes the validity of the compatibility conditions is applied. However, when studying optimality this small transient period can destroy the optimality of the control. Consequently, in order to obtain a global inverse optimality result (even in the linear case) requires an important change: instead of studying classical solutions we have to study weak solutions. The study of weak solutions allows the use of discontinuous inputs which is also required by other aspects of optimal control theory (e.g., the principle of optimality). On the other hand, the use of weak solutions comes with a price: the constructed CLFs may not be differentiable for a weak solution and the loss of regularity may also imply the loss of some stability properties (in strong norms). Therefore, the regularity properties of the weak solution must be selected very carefully: the weak solution must allow discontinuous inputs and should not involve input-dependent compatibility conditions but must be regular enough in order to retain differentiability of the CLF and some crucial stability properties. Here it should also be noted what has been observed in \cite{13, 17}: a CLF that provides estimates in weak norms and allows the feedback design on larger spaces is a feature of the CLF methodology in PDEs. Stronger stability estimates can be obtained by using additional functionals (see \cite{13, 17}).

The present paper provides the careful study of inverse optimality for the linearized model of a tank-liquid system. The model was presented in \cite{20} where a Lyapunov-based feedback design methodology was presented. The model arises from the linearization of the nonlinear model of a moving tank containing a viscous liquid that satisfies the viscous Saint-Venant equations with surface tension and wall-friction (see \cite{16, 17, 18, 21} for results dealing with the nonlinear model). The feedback design is performed in \cite{20} by using a CLF with high-order spatial derivatives and the obtained stabilization result deals with classical solutions obtained by applying the Hille-Yosida theorem. As noted above, the results obtained in \cite{20} are unusable for an inverse optimality study. Consequently, in the present work we perform the following steps:
\begin{enumerate}
\item  We define a notion of weak solution of the problem (see Definition \ref{definition1} below). The notion allows the study of discontinuous (measurable) inputs and -as expected- coincides with the classical solution when the weak solution is sufficiently regular (see Proposition \ref{proposition1} below).

\item  We show that the weak solution is unique and exists for all initial conditions and inputs (see Theorem \ref{theorem1} below). The existence/uniqueness result is obtained by applying the Galerkin method. Moreover, we show that the weak solution is regular enough to guarantee the validity of a specific equation for an appropriate Lyapunov functional (see equation \eqref{GrindEQ__13_} below). As remarked above, the Lyapunov functional that we use does not involve high-order spatial derivatives.

\item  We exploit the equation for the Lyapunov functional and we construct stabilizing feedback laws which are optimal for meaningful cost functionals (see Theorem \ref{theorem2} and Theorem \ref{theorem3} below). By obtaining accurate estimates, we also manage to show the "infinite gain margin" property and "gain-reduction margin" property, which are typical for inverse optimality problems.

This is the final step of the inverse optimality study that we perform. But in the present paper we also do more:

\item  We employ a methodology similar to the methodology for proving Input-to-State Stability for parabolic PDEs (see for example \cite{14}) and we show that the optimal feedback laws guarantee stronger stability estimates (similar to the estimates provided in \cite{20} for classical solutions; see Theorem \ref{theorem4} below). 
\end{enumerate}

To our knowledge, this is the first time that such a complete inverse optimality study has been performed for a control system described by PDEs with boundary control.

The structure of the paper is as follows. In Section 2 we present all main results of the paper. The proofs of the main results are all provided in Section 3. Finally, Section 4 is devoted to the presentation of the concluding remarks of the paper.

\noindent \textbf{Notation}
Throughout this paper, we adopt the following notation.  
\begin{itemize}
\item[$*$] $\mathbb{R}_{+} =[0,+\infty )$ denotes the set of non-negative real numbers. 
\item[$*$] Let $X$ be a given normed linear space with norm $\left\| \cdot \right\| _{X} $, let $B\subseteq X$ be a subset of $X$ and let $I\subseteq {\mathbb R}$ be a non-empty interval. By $C^{0} \left(I;B\right)$ we denote the class of continuous mappings $f:I\to B$, i.e., the class of mappings $f:I\to B$ for which the following property holds: for every $t_{0} \in I$ and for every $\varepsilon >0$ there exists $\delta >0$ such that $\left\| f(t)-f(t_{0} )\right\| _{X} <\varepsilon $ for all $t\in I\cap (t_{0} -\delta ,t_{0} +\delta )$. By $C_{c}^{0} \left(I;X\right)$ we denote the class of continuous mappings $f:I\to X$ with compact support. By $C_{c}^{1} \left(I\right)$, where $I\subseteq {\mathbb R}$ is an open interval, we denote the class of continuously differentiable functions $f:I\to {\mathbb R}$ with compact support. By $C^{1} \left(I;B\right)$ we denote the class of continuous mappings $f\in C^{0} \left(I;B\right)$ for which the following property holds: there exists $g\in C^{0} \left(I;X\right)$ such for every $t_{0} \in I$ and for every $\varepsilon >0$ there exists $\delta >0$ with $\left\| \frac{f(t)-f(t_{0} )}{t-t_{0} } -g(t_{0} )\right\| _{X} <\varepsilon $ for all $t\in I\cap \left((t_{0} -\delta ,t_{0} )\cup (t_{0} ,t_{0} +\delta )\right)$. In this case we write $f_{t} =g$. A mapping $f:I\to B$ is of class $C^{2} \left(I;B\right)$ if $f\in C^{1} \left(I;B\right)$ and $f_{t} \in C^{1} \left(I;X\right)$. 
\item[$*$] Let $S\subseteq {\mathbb R}^{n} $ be an open set and let $A\subseteq {\mathbb R}^{n} $ be a set that satisfies $S\subseteq A\subseteq cl(S)$. By $C^{0} (A\; ;\; \Omega )$, we denote the class of continuous functions on $A$, which take values in $\Omega \subseteq {\mathbb R}^{m} $. By $C^{k} (A\; ;\; \Omega )$, where $k\ge 1$ is an integer, we denote the class of functions on $A\subseteq \Re ^{n} $, which take values in $\Omega \subseteq \Re ^{m} $ and have continuous derivatives of order $k$. In other words, the functions of class $C^{k} (A;\Omega )$ are the functions which have continuous derivatives of order $k$ in $S=int(A)$ that can be continued continuously to all points in $\partial S\cap A$.  When $\Omega ={\mathbb R}$ then we write $C^{0} (A\; )$ or $C^{k} (A\; )$. When $I\subseteq {\mathbb R}$ is an interval and $\eta \in C^{1} (I)$ is a function of a single variable, $\eta '(x)$ denotes the derivative with respect to $x\in I$.
\item[$*$] Let  $I\subseteq {\mathbb R}$ be an interval and let $a<b$ be given constants. Let the function $u:I\times (a,b)\to {\mathbb R}$ be given. We use the notation $u[t]$ to denote the profile at certain $t\in I$, i.e., $(u[t])(x)=u(t,x)$ for all $x\in (a,b)$.
\item[$*$] Let $\Omega \subseteq {\mathbb R}^{n} $ be a given open set given. For $p\in [1,+\infty )$, $L^{p} (\Omega )$ is the set of equivalence classes of Lebesgue measurable functions $u:\Omega \to {\mathbb R}$ with $\left\| u\right\| _{p} :=\left(\int _{\Omega }\left|u(x)\right|^{p} dx \right)^{1/p} <+\infty $. $L^{\infty } (\Omega )$ is the set of equivalence classes of Lebesgue measurable functions $u:\Omega \to {\mathbb R}$ with $\left\| u\right\| _{\infty } :={\mathop{\text{ess}\sup }\limits_{x\in \Omega }} \left(\left|u(x)\right|\right)<+\infty $. The scalar product in $L^{2} (0,1)$ is denoted by $\left(\bullet ,\bullet \right)$. For a set $A\subseteq {\mathbb R}^{n} $, $\chi _{A} $ denotes the characteristic function of $A\subseteq {\mathbb R}^{n} $, i.e., $\chi _{A} (x)=1$ if $x\in A$ and $\chi _{A} (x)=0$ if $x\notin A$. $L_{loc}^{2} \left({\mathbb R}_{+} \right)$ denotes the class of Lebesgue measurable functions $f:\left(0,+\infty \right)\to {\mathbb R}$ with $\int _{0}^{T}\left|f(t)\right|^{2} dt <+\infty $ for all $T>0$.
\item[$*$] Let $X$ be a given Banach space with norm $\left\| \cdot \right\| _{X} $ and let $I\subseteq {\mathbb R}$ be a given open interval. A function $f:I\to X$ is measurable if there exists a set $E\subset I$ of measure zero and a sequence $\left\{f_{n} \in C_{c}^{0} \left(I;X\right)\, :\, n\ge 1\right\}$ such that ${\mathop{\lim }\limits_{n\to +\infty }} \left(f_{n} (t)\right)=f(t)$ for all $t\in I\backslash E$. For $p\in [1,+\infty )$, $L^{p} (I;X)$ is the set of equivalence classes of measurable functions $f:I\to X$ with $\left\| f\right\| _{p} :=\left(\int _{I}\left\| f(t)\right\| _{X}^{p} dx \right)^{1/p} <+\infty $. $L^{\infty } (I;X)$ is the set of equivalence classes of measurable functions $f:I\to X$ with $\left\| f\right\| _{\infty } :={\mathop{ess\sup }\limits_{t\in I}} \left(\left\| f(t)\right\| _{X} \right)<+\infty $. Weak derivatives of functions $f:I\to X$ are defined and used in \cite{4,8}. 
\item[$*$] Let $a<b$ be given constants and $k\ge 1$ be an integer. For $p\in [1,+\infty ]$, $W^{k,p} (a,b)$ denotes the Sobolev space of functions in $L^{p} (a,b)$ with weak derivatives up to order $k$, all in $L^{p} (a,b)$. We set $H^{k} (a,b)=W^{k,2} (a,b)$. The closure of $C_{c}^{1} \left((a,b)\right)$ in $H^{1} (a,b)$ is denoted by $H_{0}^{1} (a,b)$.
\end{itemize}
  
\section{Main Results}

In this paper we study the following linear PDE-ODE model for $t\ge 0$, $x\in \left(0,1\right)$
\begin{gather}
\dot{\xi }=w\quad ,\quad \dot{w}=-f \label{GrindEQ__1_} \\
\varphi _{tt} =\varphi _{xx} -\sigma \varphi _{xxxx} +\mu \varphi _{txx} -\kappa \varphi _{t} \label{GrindEQ__2_}  \\
\varphi _{x} (t,0)=\varphi _{x} (t,1)=0 \label{GrindEQ__3_}  \\
\varphi _{xxx} (t,0)=\varphi _{xxx} (t,1)=-\sigma ^{-1} f(t) \label{GrindEQ__4_} \\
\int _{0}^{1}\varphi (t,x)dx =\int _{0}^{1}\varphi _{t} (t,x)dx =0. \label{GrindEQ__5_}                                                     
\end{gather}
Model \eqref{GrindEQ__1_}-\eqref{GrindEQ__5_} is the dimensionless model of the linearization of a tank-liquid system (waterbed model) where the liquid satisfies the viscous Saint-Venant model with surface tension and wall friction (see \cite{20}).

We next provide a notion of weak solution for the problem \eqref{GrindEQ__1_}-\eqref{GrindEQ__5_}. The notion allows the study of discontinuous (measurable) inputs.

We define the following spaces:
\begin{equation} \label{GrindEQ__6_} 
\bar{S}=\left\{\, \varphi \in H^{2} (0,1)\, :\, \varphi '(0)=\varphi '(1)=\int _{0}^{1}\varphi (x)dx =0\; \right\} 
\end{equation} 
\begin{equation} \label{GrindEQ__7_} 
S=\left\{\, \varphi \in L^{2} (0,1)\, :\, \int _{0}^{1}\varphi (x)dx =0\; \right\} 
\end{equation}
Let $\mho :S\to H_{0}^{1} (0,1)$ be the linear operator defined by the following formula for all $\varphi \in S$:
\begin{equation}\label{GrindEQ__8_}
\left(\mho \varphi \right)(x)=\int _{0}^{x}\varphi (s)ds, \textrm{for all } x\in [0,1]                                              
\end{equation}
\begin{definition} \label{definition1}
Let $T>0$, $\varphi _{0} \in \bar{S}$, $\bar{\varphi }_{0} \in S$, $\left(\xi _{0} ,w_{0} \right)\in {\mathbb R}^{2} $, $f\in L^{2} (0,T)$ be given. We say that the absolutely continuous functions $\xi ,w:\left[0,T\right]\to {\mathbb R}$ and $\varphi :\left[0,T\right]\to \bar{S}$ constitute a weak solution on $[0,T]$ of the initial-boundary value problem (1)-\eqref{GrindEQ__5_} with initial condition
\begin{equation} \label{GrindEQ__9_} 
\varphi [0]=\varphi _{0} ,\varphi _{t} [0]=\bar{\varphi }_{0} , \left(\xi (0),w(0)\right)=\left(\xi _{0} ,w_{0} \right) 
\end{equation}
if \eqref{GrindEQ__1_} holds for $t\in \left(0,T\right)$ a.e., $\varphi \in L^{2} \left((0,T);H^{3} (0,1)\cap \bar{S}\right)$, \eqref{GrindEQ__9_} holds and there exists a function $u\in C^{0} \left([0,T];H_{0}^{1} (0,1)\right)\cap L^{2} \left((0,T);H^{2} (0,1)\right)$ with $u_{t} \in L^{2} \left((0,T);L^{2} (0,1)\right)$ such that
\begin{gather}
\varphi _{t} [t]=u_{x} [t], \textrm{for } t\in [0,T] \label{GrindEQ__10_} \\
u_{t} [t]=\varphi _{x} [t]-\sigma \varphi _{xxx} [t]+\mu u_{xx} [t]-\kappa u[t]-f(t)\chi _{[0,1]}, \textrm{for } t\in (0,T) \textrm{ a.e.} \label{GrindEQ__11_}                         
\end{gather}
\end{definition}
\begin{definition} \label{definition2}
Let $\varphi _{0} \in \bar{S}$, $\bar{\varphi }_{0} \in S$, $\left(\xi _{0} ,w_{0} \right)\in {\mathbb R}^{2} $, $f\in L_{loc}^{2} \left({\mathbb R}_{+} \right)$ be given. We say that the functions $\xi ,w:{\mathbb R}_{+} \to {\mathbb R}$ and $\varphi :{\mathbb R}_{+} \to \bar{S}$ constitute a weak solution on ${\mathbb R}_{+} $ of the initial-boundary value problem (1)-\eqref{GrindEQ__5_} with initial condition (9) if for every $T>0$, the functions $\xi ,w$ and $\varphi $ constitute a weak solution on $[0,T]$ of the initial-boundary value problem (1)-\eqref{GrindEQ__5_} with initial condition \eqref{GrindEQ__9_}.
\end{definition}
\begin{remark}
\begin{itemize}
\item[(i)] It should be noticed that equation \eqref{GrindEQ__10_}, definition \eqref{GrindEQ__7_} and the fact that $u\in C^{0} \left([0,T];H_{0}^{1} (0,1)\right)\cap L^{2} \left((0,T);H^{2} (0,1)\right)$ implies that $\varphi _{t} \in C^{0} \left([0,T];S\right)\cap L^{2} \left((0,T);H^{1} (0,1)\right)$. Consequently, Theorem 4 on pages 287-288 in \cite{8} implies that $\varphi $ is of class $C^{0} \left(\left[0,T\right];\bar{S}\right)$. Thus, the equalities $\varphi [0]=\varphi _{0} $, $\varphi _{t} [0]=\bar{\varphi }_{0} $ make sense.
\item[(ii)] The function $u\in C^{0} \left([0,T];H_{0}^{1} (0,1)\right)\cap L^{2} \left((0,T);H^{2} (0,1)\right)$ that satisfies \eqref{GrindEQ__10_}, \eqref{GrindEQ__11_} is strongly related to the velocity of the fluid of the original fluid model (viscous Saint-Venant model) that gives rise to model \eqref{GrindEQ__2_}-\eqref{GrindEQ__5_} (see \cite{20}).
\item[(iii)] The reader may wonder what is the inspiration behind Definition \ref{definition1} and Definition \ref{definition2}. Indeed, as said above model \eqref{GrindEQ__1_}, \eqref{GrindEQ__10_}, \eqref{GrindEQ__11_} arises as the linearization of the Saint-Venant model of an incompressible viscous fluid in a moving tank with surface tension and wall friction (see \cite{20}). In order to obtain model \eqref{GrindEQ__1_}-\eqref{GrindEQ__5_} from model \eqref{GrindEQ__1_}, \eqref{GrindEQ__10_}, \eqref{GrindEQ__11_} one needs to assume sufficient regularity and perform additional differentiations with respect to the time variable $t$ and the spatial variable $x$ (and thus eliminate the velocity from the model). 
Here we avoid the extra differentiations and consider strong solutions (in the sense of \cite{X}) of model \eqref{GrindEQ__1_}, \eqref{GrindEQ__10_}, \eqref{GrindEQ__11_}. Therefore, the notion of the weak solution provided by Definition 2.1 and Definition 2.2 is entirely motivated by the physics of the underlying problem.    
\end{itemize}
\end{remark}

The notion of weak solution which is provided by Definition \ref{definition1} coincides with the classical solution (which is guaranteed to exist by virtue of Theorem 1 in \cite{20} under certain compatibility and regularity conditions) when the weak solution is sufficiently regular. This is shown by the following proposition.

\begin{proposition} \label{proposition1}
Let $T>0$, $\varphi _{0} \in \bar{S}$, $\bar{\varphi }_{0} \in S$, $\left(\xi _{0} ,w_{0} \right)\in {\mathbb R}^{2} $, $f\in L^{2} (0,T)$ be given. Let $\xi ,w:\left[0,T\right]\to {\mathbb R}$, $\varphi \in L^{2} \left((0,T);H^{3} (0,1)\cap \bar{S}\right)$ be a weak solution on $[0,T]$ of the initial-boundary value problem \eqref{GrindEQ__1_}-\eqref{GrindEQ__5_} with initial condition \eqref{GrindEQ__9_}. Suppose that $\varphi \in C^{0} \left(\left[0,T\right];\bar{S}\cap H^{4} (0,1)\right)\cap C^{1} \left(\left[0,T\right];\bar{S}\right)\cap C^{2} \left(\left[0,T\right];L^{2} (0,1)\right)$. Then equations \eqref{GrindEQ__1_}-\eqref{GrindEQ__5_} hold for all $t\in [0,T]$.
\end{proposition}

We next show that for all initial conditions and inputs the weak solution exists and is unique. 

Define the following family of functionals for all $\varphi \in \bar{S}$, $u\in H_{0}^{1} \left(0,1\right)$ and $r\in {\mathbb R}$:
\begin{equation} \label{GrindEQ__12_} 
W\left(\varphi ,u\right):=\frac{1}{2} \left\| u\right\| _{2}^{2} +\frac{1+r}{2} \left\| \varphi \right\| _{2}^{2} +\frac{\sigma (1+r)}{2} \left\| \varphi '\right\| _{2}^{2} +\frac{r}{2} \left\| u-\mu \varphi '+\kappa \mho \varphi \right\| _{2}^{2}  
\end{equation}

\begin{theorem} \label{theorem1}
Let $T>0$, $\varphi _{0} \in \bar{S}$, $\bar{\varphi }_{0} \in S$, $\left(\xi _{0} ,w_{0} \right)\in {\mathbb R}^{2} $, $f\in L^{2} (0,T)$ be given. Then there exists a unique weak solution on $[0,T]$ of the initial-boundary value problem \eqref{GrindEQ__1_}-\eqref{GrindEQ__5_} with initial condition \eqref{GrindEQ__9_}. Moreover, for every $r\in {\mathbb R}$ the following equation holds for all $t\in [0,T]$:
\begin{align} \label{GrindEQ__13_} 
&W\left(\varphi [t],u[t]\right) 
+\int _{0}^{t}\left(r\left(\mu +\kappa \sigma \right)\left\| \varphi _{x} [s]\right\| _{2}^{2} +r\mu \sigma \left\| \varphi _{xx} [s]\right\| _{2}^{2} \right ) ds \nonumber \\
&+\int _{0}^{t}  \left( \kappa r\left\| \varphi [s]\right\| _{2}^{2} +\kappa \left\| u[s]\right\| _{2}^{2} +\mu \left\| u_{x} [s]\right\| _{2}^{2} \right)ds \nonumber  \\
=&W\left(\varphi _{0} ,\mho \bar{\varphi }_{0} \right) \nonumber \\
&+r\int _{0}^{t}f(s)\left(\int _{0}^{1}\left(\kappa x\varphi (s,x)-\frac{r+1}{r} u(s,x)+\mu \varphi _{x} (s,x)\right)dx \right)ds  
\end{align}
where 
\begin{equation*}
u\in C^{0} \left([0,T];H_{0}^{1} (0,1)\right)\cap L^{2} \left((0,T);H^{2} (0,1)\right) \textrm{with } u_{t} \in L^{2} \left((0,T);L^{2} (0,1)\right)
\end{equation*} 
is the function that satisfies \eqref{GrindEQ__10_}, \eqref{GrindEQ__11_}.
\end{theorem}
\begin{remark}
The regularity properties of the weak solution allow the validity of equation \eqref{GrindEQ__13_}. Equation \eqref{GrindEQ__13_} is crucial for the construction of a stabilizing feedback and the inverse optimality study.
\end{remark}
We next define the instantaneous cost functional which is used in our inverse optimality study and we show that it is indeed a meaningful cost. 
Let $\gamma,r,k,Q>0$ be given constants. We define for all $\varphi \in \bar{S}$, $u\in H_{0}^{1} (0,1)$, $\left(\xi ,w\right)\in {\mathbb R}^{2} $ the following functional    
\begin{align} \label{GrindEQ__40_} 
q\left(\xi ,w,\varphi ,u\right):=&k\xi ^{2} +\frac{\gamma -2 }{2k } \left(w+k\xi \right)^{2} +Q\kappa \left\| u\right\| _{2}^{2} +Q\mu \left\| u'\right\| _{2}^{2}  \nonumber \\ 
&+Q\kappa r\left\| \varphi \right\| _{2}^{2} +Qr\left(\sigma \kappa +\mu \right)\left\| \varphi '\right\| _{2}^{2} +Q\mu \sigma r\left\| \varphi ''\right\| _{2}^{2}  \nonumber \\ 
& +\frac{\gamma k^{3}  r^{2} }{2} Q^{2} \left(\int _{0}^{1}\left(\kappa x\varphi (x)-\frac{r+1}{r} u(x)+\mu \varphi '(x)\right)dx \right)^{2}  \nonumber \\ 
&-\gamma k  Qr \left(w+k\xi \right)\int _{0}^{1}\left(\kappa x\varphi (x)-\frac{r+1}{r} u(x)+\mu \varphi '(x)\right)dx  
\end{align} 
For the constants $r,k,Q>0$ we make the following assumption.

\quad

\noindent \textbf{(A)} \textit{The constants $r,k,Q>0$ satisfy the following inequality}
\begin{align} \label{GrindEQ__41_} 
k^{3}Q<\frac{1+ \sigma \pi ^{2}}{\mu r \displaystyle{\left( \left(1+\frac{\kappa}{\pi \mu \sqrt{3}}\right)^{2}+\frac{(r+1)^{2}(1+\sigma \pi ^ {2})}{r \mu (\kappa + \mu \pi ^{2})} \right)}} 
\end{align}
Assumption (A) allows us to define the critical gain
\begin{align} \label{GrindEQ__A1_} 
\gamma ^ {*}:=\frac{2}{\displaystyle{1-\frac{\mu r k^{3}Q}{1+ \sigma \pi ^{2}}\left( \left(1+\frac{\kappa}{\pi \mu \sqrt{3}}\right)^{2}+\frac{(r+1)^{2}(1+\sigma \pi ^ {2})}{r \mu (\kappa + \mu \pi ^{2})} \right)}}
\end{align}
which is used in what follows.

The fact that the instantaneous cost defined by \eqref{GrindEQ__40_} is a quadratic, coercive cost under assumption (A) for sufficiently large $\gamma$ is guaranteed by the following lemma.
\begin{lemma} \label{lemma1}
Suppose that assumption (A) holds and suppose that $\gamma >\gamma^{*}$, where $\gamma^{*}$ is defined by \eqref{GrindEQ__A1_}. Then there exists a constant $A>0$ such that
\begin{equation} \label{GrindEQ__42_}
q\left(\xi ,w,\varphi ,u\right)\ge A\left(\xi ^{2} +w^{2} +\left\| u\right\| _{H^{1} (0,1)}^{2} +\left\| \varphi \right\| _{H^{2} (0,1)}^{2} \right), 
\end{equation}
for all $\varphi \in \bar{S}$, $u\in H_{0}^{1} (0,1)$, $\left(\xi ,w\right)\in {\mathbb R}^{2}$.
\end{lemma}

We next exploit equation \eqref{GrindEQ__13_} and we construct a Control Lyapunov Functional for system \eqref{GrindEQ__1_}-\eqref{GrindEQ__5_} which allows the construction of stabilizing feedback laws of LgV type. 

Define for all $\gamma,r,k,Q>0$, $\varphi \in \bar{S}$, $\bar{\varphi }\in S$, $\left(\xi ,w\right)\in {\mathbb R}^{2} $ the following linear functional
\begin{align}  
P\left(\xi ,w,\varphi ,\bar{\varphi }\right) 
:=&\gamma k^{3} \Bigg( k^{-2}w+k^{-1} \xi  \nonumber  \\
& - rQ  \int _{0}^{1}\left(\kappa x\varphi (x)+\frac{r+1}{r} x\bar{\varphi }(x)+\mu \varphi '(x)\right)dx \Bigg ) \label{GrindEQ__56_} 
\end{align} 
and define for all $\varphi \in \bar{S}$, $u\in H_{0}^{1} (0,1)$, $\left(\xi ,w\right)\in {\mathbb R}^{2} $ the following functional
\begin{equation} \label{GrindEQ__57_} 
V\left(\xi ,w,\varphi ,u\right):=\frac{1}{2} \xi ^{2} +\frac{1}{2k^{2} } \left(w+k\xi \right)^{2} +QW\left(\varphi ,u\right) 
\end{equation} 
where $W\left(\varphi ,u\right)$ is defined by \eqref{GrindEQ__12_}. We also provide the following lemma.
\begin{lemma} \label{lemma2}
Suppose that assumption (A) holds and that $\gamma> \frac{\gamma^{*}}{2}$
where $\gamma ^{*}$ is defined by \eqref{GrindEQ__A1_}. Then, there exists a constant $\overline{A}>0$ such that
\begin{equation} \label{GrindEQ__B1_}
q\left(\xi ,w,\varphi ,u\right) + \frac{P^{2}\left(\xi,w,\varphi,u '\right)}{2\gamma k^{3}}
\ge \bar{A}\left(\xi ^{2} +w^{2} +\left\| u\right\| _{H^{1} (0,1)}^{2} +\left\| \varphi \right\| _{H^{2} (0,1)}^{2} \right), 
\end{equation}
for all $\varphi \in \bar{S}$, $u\in H_{0}^{1} (0,1)$, $\left(\xi ,w\right)\in {\mathbb R}^{2}$, where $P\left(\xi,w,\varphi,\bar{\varphi} \right) $ is the functional defined by \eqref{GrindEQ__56_}.
\end{lemma}
The proof of Lemma \ref{lemma2} is a direct consequence of Lemma \ref{lemma1} because the functional
\begin{align*}
q\left(\xi ,w,\varphi ,u\right) + \frac{P^{2}\left(\xi,w,\varphi,u' \right)}{2\gamma k^{3}}
\end{align*}
is equal to the right-hand side of \eqref{GrindEQ__40_} with $\gamma$ replaced by $2\gamma$.
\begin{theorem} \label{theorem2}
Let $r,k,Q>0$ be constants for which assumption (A) holds and $\gamma>\frac{\gamma ^{*}}{2}$, where $\gamma^{*}$ is defined by \eqref{GrindEQ__A1_}. Let $T>0$, $\varphi _{0} \in \bar{S}$, $\bar{\varphi }_{0} \in S$, $\left(\xi _{0} ,w_{0} \right)\in {\mathbb R}^{2} $ be given. Then there exists a unique weak solution on $[0,T]$ of the initial-boundary value problem \eqref{GrindEQ__1_}-\eqref{GrindEQ__5_} with initial condition \eqref{GrindEQ__9_} which satisfies for $t\in [0,T]$
\begin{equation} \label{GrindEQ__58_} 
f(t)=P\left(\xi (t),w(t),\varphi [t],\varphi _{t} [t]\right) 
\end{equation} 
where $P\left(\xi,w,\varphi,\bar{\varphi} \right) $ is the functional defined by \eqref{GrindEQ__56_}. Furthermore, there exists a constant $c>0$ that depends only on $\gamma,r,k,Q,\mu ,\sigma >0$ and $\kappa \ge 0$ (i.e., is independent of $T>0$, $\varphi _{0} \in \bar{S}$, $\bar{\varphi }_{0} \in S$, $\left(\xi _{0} ,w_{0} \right)\in {\mathbb R}^{2} $) such that the following inequality holds for all $t\in [0,T]$    
\begin{equation} \label{GrindEQ__59_} 
V\left(\xi (t),w(t),\varphi [t],u[t]\right)\le \exp \left(-ct\right)V\left(\xi _{0} ,w_{0} ,\varphi _{0} ,\mho \bar{\varphi }_{0} \right) 
\end{equation} 
where 
$$u\in C^{0} \left([0,T];H_{0}^{1} (0,1)\right)\cap L^{2} \left((0,T);H^{2} (0,1)\right) \textrm{with } u_{t} \in L^{2} \left((0,T);L^{2} (0,1)\right)$$ 
is the function that satisfies \eqref{GrindEQ__10_}, \eqref{GrindEQ__11_}. Finally, the following equation holds for all $t\in [0,T]$
\begin{equation} \label{GrindEQ__60_} 
\begin{array}{l}
\displaystyle{V\left(\xi (t),w(t),\varphi [t],u[t]\right)+\int _{0}^{t}\left(q\left(\xi (s),w(s),\varphi [s],u[s]\right)+\frac{f^{2} (s)}{2\gamma k^{3} } \right)ds }\\ 
=V\left(\xi _{0} ,w_{0} ,\varphi _{0} ,\mho \bar{\varphi }_{0} \right)
\end{array} 
\end{equation} 
where $q\left(\xi ,w,\varphi ,u\right)$ is the functional defined by \eqref{GrindEQ__40_}.
\end{theorem}  
\begin{remark}
\begin{itemize}
\item[(i)] The stability estimate \eqref{GrindEQ__59_} guarantees global exponential stabilization in the norm $\left(\left|\xi \right|^{2} +\left|w\right|^{2} +\left\| \mho \varphi _{t} \right\| _{2}^{2} +\left\| \varphi \right\| _{H^{1} (0,1)}^{2} \right)^{1/2} $ (recall (8), \eqref{GrindEQ__10_}, \eqref{GrindEQ__12_} and \eqref{GrindEQ__57_}). In other words, the stability estimate \eqref{GrindEQ__59_} in conjunction with \eqref{GrindEQ__12_} and \eqref{GrindEQ__57_} guarantees the existence of a constant $M>0$ such that for every  $\varphi _{0} \in \bar{S}$, $\bar{\varphi }_{0} \in S$, $\left(\xi _{0} ,w_{0} \right)\in {\mathbb R}^{2} $ the corresponding weak solution on ${\mathbb R}_{+} $ of the initial-boundary value problem \eqref{GrindEQ__1_}-\eqref{GrindEQ__5_} with initial condition \eqref{GrindEQ__9_} which satisfies \eqref{GrindEQ__58_} for $t\ge 0$ also satisfies the following estimate for all $t\ge 0$:
\begin{equation} \label{GrindEQ__61_} 
\begin{array}{l} {\left|\xi (t)\right|^{2} +\left|w(t)\right|^{2} +\left\| \mho \varphi _{t} [t]\right\| _{2}^{2} +\left\| \varphi [t]\right\| _{H^{1} (0,1)}^{2} } \\ {\le M\exp (-ct)\left(\left|\xi _{0} \right|^{2} +\left|w_{0} \right|^{2} +\left\| \mho \bar{\varphi }_{0} \right\| _{2}^{2} +\left\| \varphi _{0} \right\| _{H^{1} (0,1)}^{2} \right)} \end{array} 
\end{equation} 
The reader should notice that this is a weaker stabilization result than the stabilization result of Theorem 4 in \cite{20} which guarantees global exponential stabilization in the norm $\left(\left|\xi \right|^{2} +\left|w\right|^{2} +\left\| \varphi _{t} \right\| _{2}^{2} +\left\| \varphi \right\| _{H^{2} (0,1)}^{2} \right)^{1/2} $ (a stronger norm since $\left\| \varphi \right\| _{H^{1} (0,1)}^{2} \le \left\| \varphi \right\| _{H^{2} (0,1)}^{2} $ and since definition \eqref{GrindEQ__8_} implies the inequality $\left\| \mho \varphi _{t} \right\| _{2}^{2} \le \left\| \varphi _{t} \right\| _{2}^{2} $) for classical solutions of the closed-loop system \eqref{GrindEQ__1_}-\eqref{GrindEQ__5_} with \eqref{GrindEQ__58_}. This phenomenon is justified by the fact that the CLF defined by \eqref{GrindEQ__57_} involves no high-order spatial derivatives of the state. This happens because we require the feedback law to be of LgV type- a CLF with high-order spatial derivatives would require high-order spatial derivatives to appear in a LgV type feedback law. Moreover, a CLF containing terms like $\left\| \varphi _{t} \right\| _{2}^{2} $ or $\left\| \varphi \right\| _{H^{2} (0,1)}^{2} $ (see the CLF used in \cite{20}) is not necessarily differentiable with respect to time for a weak solution.
\item[(ii)] The feedback law \eqref{GrindEQ__56_} is different from the feedback law proposed in \cite{20} because it contains the term $\kappa \int _{0}^{1}x\varphi (x)dx $. Therefore, the feedback law \eqref{GrindEQ__56_} depends on the friction coefficient $\kappa \ge 0$ in contrast with the feedback law proposed in \cite{20}. Thus, the feedback law proposed in \cite{20} guarantees robustness with respect to the friction coefficient. Of course, the feedback law proposed in \cite{20} is not optimal for a meaningful cost criterion. 
\item[(iii)] The proof of Theorem \ref{theorem2}, provided in the following section, shows that every weak solution of \eqref{GrindEQ__1_}-\eqref{GrindEQ__5_} with \eqref{GrindEQ__58_} satisfies for $t \in (0,T)$ a.e. the following equation
\begin{align*}  \label{C1}   
&\frac{d}{d\, t} V\left(\xi (t),w(t),\varphi [t],u[t]\right) \nonumber \\
=&\xi (t)w(t)+k^{-1} \left(w(t)+k\xi (t)\right)w(t)-Q\kappa \left\| u[t]\right\| _{2}^{2} -Q\mu \left\| u_{x} [t]\right\| _{2}^{2} \nonumber  \\ 
&-Q\kappa r\left\| \varphi [t]\right\| _{2}^{2} -Qr\left(\sigma \kappa +\mu \right)\left\| \varphi _{x} [t]\right\| _{2}^{2} -Q\mu \sigma r\left\| \varphi _{xx} [t]\right\| _{2}^{2} \nonumber \\ 
&-f(t)\left(k^{-2} \left(w(t)+k\xi (t)\right)-Qr\kappa \int _{0}^{1} x \varphi(t,x) dx \right) \nonumber \\
&+f(t)\left(Q(r+1)\int _{0}^{1}x \varphi_{t}(t,x)dx +Qr\mu \int_{0}^{1} \varphi_{x}(t,x) dx \right)
\end{align*}
Therefore, the feedback law \eqref{GrindEQ__58_} with $P\left(\xi,w,\varphi,\bar{\varphi} \right) $ defined by \eqref{GrindEQ__56_}  is a LgV type feedback.  
\end{itemize}
\end{remark}

The inverse optimality study is finished by showing that the stabilizing feedback law \eqref{GrindEQ__58_} is optimal for a meaningful cost functional. By the term ``meaningful cost functional" for a linear control system we mean the sum of two quadratic, coercive functionals: one for the state and one for the control input. For every $\varphi _{0} \in \bar{S}$, $\bar{\varphi }_{0} \in S$, $\left(\xi _{0} ,w_{0} \right)\in {\mathbb R}^{2} $, $f\in L_{loc}^{2} \left({\mathbb R}_{+} \right)$ we define the functional 
\begin{equation} \label{GrindEQ__91_} 
J\left(\xi _{0} ,w_{0} ,\varphi _{0} ,\bar{\varphi }_{0} ,f\right)=\int _{0}^{+\infty }\left(q\left(\xi (s),w(s),\varphi [s],u[s]\right)+\frac{f^{2} (s)}{2\gamma k^{3} } \right)ds  
\end{equation} 
where $\xi ,w:{\mathbb R}_{+} \to {\mathbb R}$ and $\varphi :{\mathbb R}_{+} \to \bar{S}$ is the unique weak solution on ${\mathbb R}_{+} $ of the initial-boundary value problem \eqref{GrindEQ__1_}-\eqref{GrindEQ__5_} with initial condition \eqref{GrindEQ__9_} and $u:{\mathbb R}_{+} \to H_{0}^{1} (0,1)$ is the function that satisfies \eqref{GrindEQ__10_}, \eqref{GrindEQ__11_}. When $r,k,Q>0$ are constants for which assumption (A) holds and $\gamma>\gamma^{*}$, Lemma \ref{lemma1} shows that the functional $J$ defined by \eqref{GrindEQ__91_} is a meaningful cost functional. 
\begin{theorem} \label{theorem3}
Let $r,k,Q>0$ be constants for which assumption (A) holds and $\gamma>\gamma^{*}$, where $\gamma^{*}$ is defined by \eqref{GrindEQ__A1_}. Then for every $\varphi _{0} \in \bar{S}$, $\bar{\varphi }_{0} \in S$, $\left(\xi _{0} ,w_{0} \right)\in {\mathbb R}^{2} $ it holds that
\begin{align} \label{GrindEQ__92_} 
V\left(\xi _{0} ,w_{0} ,\varphi _{0} ,\mho \bar{\varphi }_{0} \right)&=\inf \left\{\, J\left(\xi _{0} ,w_{0} ,\varphi _{0} ,\bar{\varphi }_{0} ,f\right)\, :\, f\in L_{loc}^{2} \left({\mathbb R}_{+} \right)\, \right\} \nonumber \\ 
&=J\left(\xi _{0} ,w_{0} ,\varphi _{0} ,\bar{\varphi }_{0} ,f^{*} \right) 
\end{align} 
where 
\begin{equation} \label{GrindEQ__94_}
f^{*} (t)=P\left(\xi ^{*} (t),w^{*} (t),\varphi ^{*} [t],\varphi _{t}^{*} [t]\right), \textrm{for } t\ge 0, 
\end{equation}
$P\left(\xi,w,\varphi,\bar{\varphi} \right) $ is the functional defined by \eqref{GrindEQ__56_} and $\xi ^{*} ,w^{*} :{\mathbb R}_{+} \to {\mathbb R}$, $\varphi ^{*} :{\mathbb R}_{+} \to \bar{S}$ is the unique weak solution on ${\mathbb R}_{+} $ of the initial-boundary value problem \eqref{GrindEQ__1_}-\eqref{GrindEQ__5_} with initial condition \eqref{GrindEQ__9_} that satisfies \eqref{GrindEQ__94_}.
\end{theorem}
\begin{remark}
\begin{itemize}
\item[(i)] It is important to notice that equation \eqref{GrindEQ__92_} is different from the equation  
\begin{equation} \label{GrindEQ__95_} 
V\left(\xi _{0} ,w_{0} ,\varphi _{0} ,\mho \bar{\varphi }_{0} \right)=\inf \left\{\, J\left(\xi _{0} ,w_{0} ,\varphi _{0} ,\bar{\varphi }_{0} ,f\right)\, :\, f\in F\left(\xi _{0} ,w_{0} ,\varphi _{0} ,\bar{\varphi }_{0} \right)\, \right\} 
\end{equation} 
where $F\left(\xi _{0} ,w_{0} ,\varphi _{0} ,\bar{\varphi }_{0} \right)\subseteq L_{loc}^{2} \left({\mathbb R}_{+} \right)$ is the class of inputs for which the corresponding unique weak solution on ${\mathbb R}_{+} $ $\xi ,w:{\mathbb R}_{+} \to {\mathbb R}$ and $\varphi :{\mathbb R}_{+} \to \bar{S}$ of the initial-boundary value problem \eqref{GrindEQ__1_}-\eqref{GrindEQ__5_} with initial condition \eqref{GrindEQ__9_} satisfies ${\mathop{\lim }\limits_{t\to +\infty }} \left(V\left(\xi (t),w(t),\varphi [t],\mho \varphi _{t} [t]\right)\right)=0$ (see for example the analysis in \cite{33}). Therefore, Theorem \ref{theorem3} is a stronger result than the usual results of inverse optimality.
\item[(ii)] The fact that Theorem \ref{theorem3} holds for $\gamma > \gamma ^{*}$, while Theorem \ref{theorem2} holds for $\gamma > \frac{\gamma ^{*}}{2}$ is typical for inverse optimality problems. The fact that $\gamma$ is unrestricted shows the "infinite gain margin" property which holds for optimal feedback laws, while the difference between assumptions $\gamma>\gamma^{*}$ and $\gamma>\frac{\gamma^{*}}{2}$ shows the "gain-reduction margin" property.
\end{itemize}
\end{remark}
Setting $Q=0$, we obtain the following corollary. 
\begin{corollary}
For the tank \eqref{GrindEQ__1_} without liquid and with feedback $f=\gamma k (w + k \xi) $ and arbitrary $k>0$, exponential stability holds in $\|(\xi,w)\|$ for all $\gamma>1$  and the infinite-time integral of
$ \gamma \left(2(k \xi)^{2} +(\gamma-2) \left(w+k \xi \right)^{2}\right) + \left(\frac{f}{k}\right)^{2}$
is minimized for all $\gamma>2$.
\end{corollary}
As remarked in the Introduction in this work we go one step beyond the inverse optimality study. The following theorem shows that the optimal feedback law \eqref{GrindEQ__58_} guarantees the same stabilization result as the one obtained for classical solutions in \cite{20}. However, the stability estimates in stronger norms have to be obtained by exploiting additional functionals or a different procedure. Here we employ a methodology similar to the methodology for proving Input-to-State Stability for parabolic PDEs (see for example \cite{14}).
\begin{theorem} \label{theorem4}
Let $r,k,Q>0$ be constants for which assumption (A) holds and $\gamma>\frac{\gamma^{*}}{2}$, where $\gamma^{*}$ is defined by \eqref{GrindEQ__A1_}. Then there exist constants $G,\omega >0$ such that for every $\varphi _{0} \in \bar{S}$, $\bar{\varphi }_{0} \in S$, $\left(\xi _{0} ,w_{0} \right)\in {\mathbb R}^{2} $ the unique weak solution on ${\mathbb R}_{+} $ of the initial-boundary value problem \eqref{GrindEQ__1_}-\eqref{GrindEQ__5_} with initial condition \eqref{GrindEQ__9_} which satisfies \eqref{GrindEQ__58_} for $t\ge 0$ also satisfies the following estimate for all $t\ge 0$:
\begin{equation} \label{GrindEQ__103_} 
\begin{array}{l}
\left|\xi (t)\right|^{2} +\left|w(t)\right|^{2} +\left\| \varphi _{t} [t]\right\| _{2}^{2} +\left\| \varphi [t]\right\| _{H^{2} (0,1)}^{2} \\
\le G\exp \left(-\omega \, t\right)\left(\left|\xi _{0} \right|^{2} +\left|w_{0} \right|^{2} +\left\| \bar{\varphi }_{0} \right\| _{2}^{2} +\left\| \varphi _{0} \right\| _{H^{2} (0,1)}^{2} \right)
\end{array} 
\end{equation} 
\end{theorem}
\section{Proofs}
We start by providing the proof of Proposition \ref{proposition1}. 
\begin{proof}{\hspace{-0.15em}\textit{of Proposition \ref{proposition1}}}
Since $\varphi \in C^{0} \left(\left[0,T\right];\bar{S}\cap H^{4} (0,1)\right)\cap C^{1} \left(\left[0,T\right];\bar{S}\right)\cap C^{2} \left(\left[0,T\right];L^{2} (0,1)\right)$, it follows that $\varphi _{t} \in C^{0} \left(\left[0,T\right];\bar{S}\right)\cap C^{1} \left(\left[0,T\right];L^{2} (0,1)\right)$ and $\varphi _{tt} \in C^{0} \left(\left[0,T\right];L^{2} (0,1)\right)$. By virtue of \eqref{GrindEQ__10_} and definition \eqref{GrindEQ__8_} we conclude that $u[t]=\mho \varphi _{t} [t]$ for all $t\in [0,T]$ and therefore
 $u\in C^{0} \left(\left[0,T\right];H^{3} (0,1)\right.$ $\left. \cap H_{0}^{1} (0,1)\right) \cap C^{1} \left(\left[0,T\right];H^{1} (0,1)\right)$
  with $u_{xx} (t,0)=u_{xx} (t,1)=0$ for all $t\in [0,T]$. Moreover, \eqref{GrindEQ__11_} implies that $f\chi _{[0,1]} \in C^{0} \left(\left[0,T\right];L^{2} (0,1)\right)$ which gives that $f\in C^{0} \left(\left[0,T\right]\right)$. Thus both \eqref{GrindEQ__10_} and \eqref{GrindEQ__11_} hold for all $t\in [0,T]$. This fact implies that $\varphi _{tt} [t]=u_{tx} [t]$, $\varphi _{tx} [t]=u_{xx} [t]$, $\varphi _{txx} [t]=u_{xxx} [t]$, $u_{tx} [t]=\varphi _{xx} [t]-\sigma \varphi _{xxxx} [t]+\mu u_{xxx} [t]-\kappa u_{x} [t]$ for all $t\in [0,T]$. Therefore, we conclude that \eqref{GrindEQ__2_} holds for all $t\in [0,T]$.

The fact that \eqref{GrindEQ__3_} and \eqref{GrindEQ__5_} hold for all $t\in [0,T]$ follows from the fact that $\varphi ,\varphi _{t} \in C^{0} \left(\left[0,T\right];\bar{S}\right)$ (recall definition \eqref{GrindEQ__6_}). Since $f\in C^{0} \left(\left[0,T\right]\right)$, it follows that $\xi ,w\in C^{1} \left(\left[0,T\right]\right)$ and that \eqref{GrindEQ__1_} holds for all $t\in [0,T]$. Since \eqref{GrindEQ__11_} holds for all $t\in [0,T]$ and since $u_{t} [t],\varphi _{x} [t],\varphi _{xxx} [t],u_{xx} [t],u[t],f(t)\chi _{[0,1]} $ are all continuous functions on $\left[0,1\right]$, it follows that $u_{t} (t,x)=\varphi _{x} (t,x)-\sigma \varphi _{xxx} (t,x)+\mu u_{xx} (t,x)-\kappa u(t,x)-f(t)$ for all $(t,x)\in \left[0,T\right]\times \left[0,1\right]$. Consequently, since $u(t,0)=u(t,1)=\varphi _{x} (t,0)=\varphi _{x} (t,1)=u_{xx} (t,0)=u_{xx} (t,1)=0$ for all $t\in [0,T]$ (which also imply the equations $u_{t} (t,0)=u_{t} (t,1)=0$), we obtain that \eqref{GrindEQ__4_} holds for all $t\in [0,T]$.    

\noindent The proof is complete. 
\qed
\end{proof}

The proof of Theorem \ref{theorem1} is based on the Galerkin method as explained in the book \cite{8}.
\begin{proof}{\hspace{-0.15em}\textit{of Theorem \ref{theorem1}}}
The proof is divided in three parts. 

\noindent \underbar{1${}^{st}$ Part: Existence.} We prove that there exists a unique weak solution on $[0,T]$ of the initial-boundary value problem \eqref{GrindEQ__1_}-\eqref{GrindEQ__5_} with initial condition \eqref{GrindEQ__9_}. 

 We construct approximate solutions of the initial-boundary value problem \eqref{GrindEQ__1_}-\eqref{GrindEQ__5_} with initial condition \eqref{GrindEQ__9_} by means of the following formula for $m=1,2,...$:
\begin{equation} \label{GrindEQ__14_}
\varphi _{m} (t,x)=\sum _{n=1}^{m}a_{n} (t)\phi _{n} (x) , \textrm{for } (t,x)\in \left[0,T\right]\times \left[0,1\right]
\end{equation}
where
\begin{equation} \label{GrindEQ__15_}
\phi _{n} (x)=\sqrt{2} \cos \left(n\pi x\right), \textrm{for } x\in \left[0,1\right]                                             
\end{equation}
and each $a_{n} $ for $n=1,...,m$ is the solution of the following initial value problem
\begin{equation} \label{GrindEQ__16_} 
\ddot{a}_{n} (t)=-\left(\mu n^{2} \pi ^{2} +\kappa \right)\dot{a}_{n} (t)-\left(\sigma n^{2} \pi ^{2} +1\right)n^{2} \pi ^{2} a_{n} (t)+\sqrt{2} f(t)\left((-1)^{n} -1\right) 
\end{equation} 
\begin{equation} \label{GrindEQ__17_} 
a_{n} (0)=\sqrt{2} \int _{0}^{1}\varphi _{0} (x)\cos \left(n\pi x\right)dx , \dot{a}_{n} (0)=\sqrt{2} \int _{0}^{1}\bar{\varphi }_{0} (x)\cos \left(n\pi x\right)dx  
\end{equation}
Multiplying \eqref{GrindEQ__16_} by $\dot{a}_{n} (t)$ we get for a.e. $t\in \left(0,T\right)$:
\begin{equation}\label{GrindEQ__18_}
\begin{array}{l}
\displaystyle{\frac{d}{d\, t} \left(\frac{1}{2} \dot{a}_{n}^{2} (t)+\frac{1}{2} \left(\sigma n^{2} \pi ^{2} +1\right)n^{2} \pi ^{2} a_{n}^{2} (t)\right) }\\
\displaystyle{=-\left(\mu n^{2} \pi ^{2} +\kappa \right)\dot{a}_{n}^{2} (t)+\sqrt{2} f(t)\left((-1)^{n} -1\right)\dot{a}_{n} (t)  }
\end{array}
\end{equation} 
Rewriting \eqref{GrindEQ__16_} in the following way
\[
\begin{array}{l}
\displaystyle{\frac{d}{d\, t} \left(\dot{a}_{n} (t)+\left(\mu n^{2} \pi ^{2} +\kappa \right)a_{n} (t)\right)} \\
=\displaystyle{-\left(\sigma n^{2} \pi ^{2} +1\right)n^{2} \pi ^{2} a_{n} (t)+\sqrt{2} f(t)\left((-1)^{n} -1\right)}
\end{array}
\] 
and multiplying the above equation by $\left(\dot{a}_{n} (t)+\left(\mu n^{2} \pi ^{2} +\kappa \right)a_{n} (t)\right)$ we get for a.e. $t\in \left(0,T\right)$:
\begin{align} \label{GrindEQ__19_} 
&\frac{d}{d\, t} \left(\frac{1}{2} \left(\dot{a}_{n} (t)+\left(\mu n^{2} \pi ^{2} +\kappa \right)a_{n} (t)\right)^{2} +\frac{1}{2} \left(\sigma n^{2} \pi ^{2} +1\right)n^{2} \pi ^{2} a_{n}^{2} (t)\right) \nonumber \\ 
=&-\left(\mu n^{2} \pi ^{2} +\kappa \right)\left(\sigma n^{2} \pi ^{2} +1\right)n^{2} \pi ^{2} a_{n}^{2} (t) \nonumber \\
&+\sqrt{2} f(t)\left((-1)^{n} -1\right)\left(\dot{a}_{n} (t)+\left(\mu n^{2} \pi ^{2} +\kappa \right)a_{n} (t)\right)  
\end{align}
Equation \eqref{GrindEQ__18_} gives for a.e. $t\in \left(0,T\right)$:
\begin{align} \label{GrindEQ__20_} 
\frac{d}{d\, t} \left(\frac{1}{2} \dot{a}_{n}^{2} (t)+\frac{1}{2} \left(\sigma n^{2} \pi ^{2} +1\right)n^{2} \pi ^{2} a_{n}^{2} (t)\right) \nonumber \\
\le -\frac{1}{2} \left(\mu n^{2} \pi ^{2} +\kappa \right)\dot{a}_{n}^{2} (t)+\frac{4f^{2} (t)}{\mu n^{2} \pi ^{2} +\kappa }  
\end{align} 
Integrating the differential inequality \eqref{GrindEQ__20_} we get for all $t\in \left[0,T\right]$:
\begin{align} \label{GrindEQ__21_} 
\dot{a}_{n}^{2} (t)+\left(\sigma n^{2} \pi ^{2} +1\right)n^{2} \pi ^{2} a_{n}^{2} (t)+\int _{0}^{t}\left(\mu n^{2} \pi ^{2} +\kappa \right)\dot{a}_{n}^{2} (s)ds \nonumber  \\ 
\le \dot{a}_{n}^{2} (0)+\left(\sigma n^{2} \pi ^{2} +1\right)n^{2} \pi ^{2} a_{n}^{2} (0)+\frac{8}{\mu n^{2} \pi ^{2} +\kappa } \int _{0}^{t}f^{2} (s)ds  
\end{align}
Equations \eqref{GrindEQ__14_}, \eqref{GrindEQ__15_} imply that
\begin{align}
&\left\| \varphi _{m,t} [t]\right\| _{2}^{2} =\sum _{n=1}^{m}\dot{a}_{n}^{2} (t), \left\| \varphi _{m,x} [t]\right\| _{2}^{2} =\sum _{n=1}^{m}n^{2} \pi ^{2} a_{n}^{2} (t), \nonumber  \\
&\left\| \varphi _{m,xx} [t]\right\| _{2}^{2} =\sum _{n=1}^{m}n^{4} \pi ^{4} a_{n}^{2} (t), \left\| \varphi _{m,tx} [t]\right\| _{2}^{2} =\sum _{n=1}^{m}n^{2} \pi ^{2} \dot{a}_{n}^{2} (t). \nonumber 
\end{align}
Consequently, we get from \eqref{GrindEQ__21_} for all $t\in \left[0,T\right]$:
\begin{align} \label{GrindEQ__22_} 
\left\| \varphi _{m,t} [t]\right\| _{2}^{2} +\left\| \varphi _{m,x} [t]\right\| _{2}^{2} +&\sigma \left\| \varphi _{m,xx} [t]\right\| _{2}^{2} + \mu \int _{0}^{t}\left\| \varphi _{m,tx} [s]\right\| _{2}^{2} ds  \nonumber \\
+\kappa \int _{0}^{t}\left\| \varphi _{m,t} [s]\right\| _{2}^{2} ds   \le & \left\| \varphi _{m,t} [0]\right\| _{2}^{2} +\left\| \varphi _{m,x} [0]\right\| _{2}^{2}  \nonumber +\sigma \left\| \varphi _{m,xx} [0]\right\| _{2}^{2} \\ 
&+\left(\sum _{n=1}^{m}\frac{8}{\mu n^{2} \pi ^{2} +\kappa }  \right)\int _{0}^{t}f^{2} (s)ds   
\end{align} 
Equation \eqref{GrindEQ__19_} gives for a.e. $t\in \left(0,T\right)$:
\begin{align} \label{GrindEQ__23_} 
&\frac{d}{d\, t} \left(\frac{1}{2} \left(\dot{a}_{n} (t)+\left(\mu n^{2} \pi ^{2} +\kappa \right)a_{n} (t)\right)^{2} +\frac{1}{2} \left(\sigma n^{2} \pi ^{2} +1\right)n^{2} \pi ^{2} a_{n}^{2} (t)\right) \nonumber \\ \le &-\frac{1}{2} \left(\mu n^{2} \pi ^{2} +\kappa \right)\left(\sigma n^{2} \pi ^{2} +1\right)n^{2} \pi ^{2} a_{n}^{2} (t)+\frac{1}{2} \left(\mu n^{2} \pi ^{2} +\kappa \right)\dot{a}_{n}^{2} (t) \nonumber \\ 
&+\left(\frac{\left(\mu n^{2} \pi ^{2} +\kappa \right)^{2} }{\left(\sigma n^{2} \pi ^{2} +1\right)n^{2} \pi ^{2} } +1\right)\frac{\left((-1)^{n} -1\right)^{2} }{\mu n^{2} \pi ^{2} +\kappa } f^{2} (t) 
\end{align} 
Integrating the differential inequality \eqref{GrindEQ__23_} we get for all $t\in \left[0,T\right]$:
\begin{align} \label{GrindEQ__24_} 
&\left(\dot{a}_{n} (t)+\left(\mu n^{2} \pi ^{2} +\kappa \right)a_{n} (t)\right)^{2} +\left(\sigma n^{2} \pi ^{2} +1\right)n^{2} \pi ^{2} a_{n}^{2} (t) \nonumber \\ 
&+\int _{0}^{t}\left(\mu n^{2} \pi ^{2} +\kappa \right)\left(\sigma n^{2} \pi ^{2} +1\right)n^{2} \pi ^{2} a_{n}^{2} (s)ds  \nonumber \\ 
\le & \left(\dot{a}_{n} (0)+\left(\mu n^{2} \pi ^{2} +\kappa \right)a_{n} (0)\right)^{2} \nonumber \\
&+\left(\sigma n^{2} \pi ^{2} +1\right)n^{2} \pi ^{2} a_{n}^{2} (0)+\int _{0}^{t}\left(\mu n^{2} \pi ^{2} +\kappa \right)\dot{a}_{n}^{2} (s)ds   \nonumber \\
&+\frac{8\left(\left(\mu +\kappa \right)^{2} +\sigma \right)}{\sigma \left(\mu n^{2} \pi ^{2} +\kappa \right)} \int _{0}^{t}f^{2} (s)ds
\end{align} 
Equations \eqref{GrindEQ__14_}, \eqref{GrindEQ__15_} imply that
\begin{align*}
&\left\| \varphi _{m,t} [t]-\mu \varphi _{m,xx} [t]+\kappa \varphi _{m} [t]\right\| _{2}^{2} =\sum _{n=1}^{m}\left(\dot{a}_{n} (t)+\left(\mu n^{2} \pi ^{2} +\kappa \right)a_{n} (t)\right)^{2}, \\
&\left\| \varphi _{m,t} [t]\right\| _{2}^{2} =\sum _{n=1}^{m}\dot{a}_{n}^{2} (t), \left\| \varphi _{m,x} [t]\right\| _{2}^{2} =\sum _{n=1}^{m}n^{2} \pi ^{2} a_{n}^{2} (t),  \\
& \left\| \varphi _{m,xx} [t]\right\| _{2}^{2} =\sum _{n=1}^{m}n^{4} \pi ^{4} a_{n}^{2} (t), 
\left\|  \varphi _{m,xxx} [t]\right\| _{2}^{2} =\sum _{n=1}^{m}n^{6} \pi ^{6} a_{n}^{2} (t), \\
&\left\| \varphi _{m,tx} [t]\right\| _{2}^{2} =\sum _{n=1}^{m}n^{2} \pi ^{2} \dot{a}_{n}^{2} (t).
\end{align*}
Consequently, we get from \eqref{GrindEQ__24_} for all $t\in \left[0,T\right]$:
\begin{align}  \label{GrindEQ__25_}
&\left\| \varphi _{m,t} [t]-\mu \varphi _{m,xx} [t]+\kappa \varphi _{m} [t]\right\| _{2}^{2} +\left\| \varphi _{m,x} [t]\right\| _{2}^{2}  \nonumber  \\ 
&+\sigma \left\| \varphi _{m,xx} [t]\right\| _{2}^{2} +\mu \sigma \int _{0}^{t}\left\| \varphi _{m,xxx} [s]\right\| _{2}^{2} ds  \nonumber \\ 
&+\left(\kappa \sigma +\mu \right)\int _{0}^{t}\left\| \varphi _{m,xx} [s]\right\| _{2}^{2} ds +\kappa \int _{0}^{t}\left\| \varphi _{m,x} [s]\right\| _{2}^{2} ds \nonumber \\
\le & \left\| \varphi _{m,t} [0]-\mu \varphi _{m,xx} [0]+\kappa \varphi _{m} [0]\right\| _{2}^{2} \nonumber \\  
&+\left\| \varphi _{m,x} [0]\right\| _{2}^{2} +\sigma \left\| \varphi _{m,xx} [0]\right\| _{2}^{2} +\mu \int _{0}^{t}\left\| \varphi _{m,tx} [s]\right\| _{2}^{2} ds +\kappa \int _{0}^{t}\left\| \varphi _{m,t} [s]\right\| _{2}^{2} ds \nonumber  \\ 
&+\left(\frac{\left(\mu +\kappa \right)^{2} }{\sigma } +1\right)\left(\sum _{n=1}^{m}\frac{8}{\mu n^{2} \pi ^{2} +\kappa }  \right)\int _{0}^{t}f^{2} (s)ds   
\end{align} 
Combining \eqref{GrindEQ__22_} and \eqref{GrindEQ__25_} we obtain for all $t\in \left[0,T\right]$:
\begin{align} \label{GrindEQ__26_} 
&\left\| \varphi _{m,t} [t]-\mu \varphi _{m,xx} [t]+\kappa \varphi _{m} [t]\right\| _{2}^{2} +2\left\| \varphi _{m,t} [t]\right\| _{2}^{2}   \nonumber  \\ 
&+3\left\| \varphi _{m,x} [t]\right\| _{2}^{2}+3\sigma \left\| \varphi _{m,xx} [t]\right\| _{2}^{2} \nonumber \\
& +\int _{0}^{t}\left(\mu \sigma \left\| \varphi _{m,xxx} [s]\right\| _{2}^{2} +\left(\kappa \sigma +\mu \right)\left\| \varphi _{m,xx} [s]\right\| _{2}^{2} +\kappa \left\| \varphi _{m,x} [s]\right\| _{2}^{2} \right)ds \nonumber \\
&+\int _{0}^{t}\left(\mu \left\| \varphi _{m,tx} [s]\right\| _{2}^{2} +\kappa \left\| \varphi _{m,t} [s]\right\| _{2}^{2} \right)ds  \nonumber \\
\le & \left\| \varphi _{m,t} [0]-\mu \varphi _{m,xx} [0]+\kappa \varphi _{m} [0]\right\| _{2}^{2} \nonumber \\
&+3\left\| \varphi _{m,x} [0]\right\| _{2}^{2} +3\sigma \left\| \varphi _{m,xx} [0]\right\| _{2}^{2} +2\left\| \varphi _{m,t} [0]\right\| _{2}^{2}  \nonumber \\  
&+\left(\frac{\left(\mu +\kappa \right)^{2} }{\sigma } +3\right)\left(\sum _{n=1}^{m}\frac{8}{\mu n^{2} \pi ^{2} +\kappa }  \right)\int _{0}^{t}f^{2} (s)ds
\end{align} 
Since $\int _{0}^{1}\varphi _{m} (t,x)dx =\int _{0}^{1}\varphi _{m,t} (t,x)dx =\varphi _{m,x} (t,0)=0$ for all $t\in \left[0,T\right]$ (recall \eqref{GrindEQ__14_}, \eqref{GrindEQ__15_}), it follows that $\left\| \varphi _{m} [t]\right\| _{2} \le \left\| \varphi _{m,x} [t]\right\| _{2} \le \left\| \varphi _{m,xx} [t]\right\| _{2} $ and $\left\| \varphi _{m,t} [t]\right\| _{2} \le \left\| \varphi _{m,tx} [t]\right\| _{2} $ for all $t\in \left[0,T\right]$. Moreover, equations \eqref{GrindEQ__14_}, \eqref{GrindEQ__17_} imply that $\left\| \varphi _{m,x} [0]\right\| _{2} \le \left\| \varphi '_{0} \right\| _{2} $, $\left\| \varphi _{m} [0]\right\| _{2} \le \left\| \varphi _{0} \right\| _{2} $, $\left\| \varphi _{m,xx} [0]\right\| _{2} \le \left\| \varphi ''_{0} \right\| _{2} $ and $\left\| \varphi _{m,t} [0]\right\| _{2} \le \left\| \bar{\varphi }_{0} \right\| _{2} $. Consequently, estimate \eqref{GrindEQ__26_} implies that there exists a constant $K>0$ (independent of $m$) such that the following inequality holds for all $t\in \left[0,T\right]$:
\begin{equation} \label{GrindEQ__27_} 
\left\| \varphi _{m,t} [t]\right\| _{2}^{2} +\left\| \varphi _{m} [t]\right\| _{H^{2} (0,1)}^{2} +\int _{0}^{t}\left\| \varphi _{m} [s]\right\| _{H^{3} (0,1)}^{2} ds +\int _{0}^{t}\left\| \varphi _{m,t} [s]\right\| _{H^{1} (0,1)}^{2} ds \le K 
\end{equation} 
We define for $t\in \left[0,T\right]$ and $x\in \left[0,1\right]$:
\begin{equation} \label{GrindEQ__28_} 
u_{m} (t,x)=\sum _{n=1}^{m}\frac{\dot{a}_{n} (t)}{n\pi } g_{n} (x)  
\end{equation} 
where
\begin{equation} \label{GrindEQ__29_}
g_{n} (x)=\sqrt{2} \sin \left(n\pi x\right), \textrm{for } x\in \left[0,1\right].                                             
\end{equation}
Equations \eqref{GrindEQ__18_}, \eqref{GrindEQ__14_}, \eqref{GrindEQ__15_}, \eqref{GrindEQ__16_}, \eqref{GrindEQ__29_} give for $t\in \left[0,T\right]$ and $x\in \left[0,1\right]$ 
\begin{equation} \label{GrindEQ__30_} 
u_{m} (t,x)=\int _{0}^{x}\varphi _{m,t} (t,s)ds  
\end{equation} 
and for a.e. $t\in \left(0,T\right)$:
\begin{align} \label{GrindEQ__31_} 
u_{m,t} [t]=&\varphi _{m,x} [t]-\kappa u_{m} [t]+\mu \varphi _{m,tx} [t] \nonumber \\ 
&-\sigma \varphi _{m,xxx} [t]+\sqrt{2} f(t)\sum _{n=1}^{m}\frac{(-1)^{n} -1}{n\pi } g_{n}   
\end{align}
Equation \eqref{GrindEQ__30_} implies that $u_{m} (t,0)=u_{m} (t,1)=0$, $\left\| u_{m} [t]\right\| _{2} \le \left\| \varphi _{m,t} [t]\right\| _{2} $, $\left\| u_{m,x} [t]\right\| _{2} \le \left\| \varphi _{m,t} [t]\right\| _{2} $ and $\left\| u_{m,xx} [t]\right\| _{2} \le \left\| \varphi _{m,tx} [t]\right\| _{2} $ for all $t\in \left[0,T\right]$. This fact in conjunction with \eqref{GrindEQ__27_} implies that there exists a constant $\bar{K}>0$ (independent of $m$) such that the following inequality holds for all $t\in \left[0,T\right]$:
\begin{equation} \label{GrindEQ__32_} 
\left\| u_{m} [t]\right\| _{H_{0}^{1} (0,1)}^{2} +\int _{0}^{t}\left\| u_{m} [s]\right\| _{H^{2} (0,1)}^{2} ds +\int _{0}^{t}\left\| u_{m,t} [s]\right\| _{2}^{2} ds \le \bar{K} 
\end{equation} 
Using the estimates 
\begin{align*}
\left|\varphi _{m} (t,x)-\varphi _{m} (t_{0} ,x)\right|&\le \int _{t_{0} }^{t}\left|\varphi _{m,t} (s,x)\right|ds  \\ 
&\le \left(t-t_{0} \right)^{1/2} \left(\int _{t_{0} }^{t}\left|\varphi _{m,t} (s,x)\right|^{2} ds \right)^{1/2},
\end{align*}
\begin{align*}
\left|u_{m} (t,x)-u_{m} (t_{0} ,x)\right|
&\le \int _{t_{0} }^{t}\left|u_{m,t} (s,x)\right|ds \\
&\le \left(t-t_{0} \right)^{1/2} \left(\int _{t_{0} }^{t}\left|u_{m,t} (s,x)\right|^{2} ds \right)^{1/2} 
\end{align*}
that hold for all $x\in \left[0,1\right]$, $t,t_{0} \in \left[0,T\right]$ with $t\ge t_{0} $ , we get the estimates  
\begin{align*}
\left\| \varphi _{m} [t]-\varphi _{m} [t_{0} ]\right\| _{2}^{2} &\le \left(t-t_{0} \right)\int _{0}^{T}\left\| \varphi _{m,t} [s]\right\| _{2}^{2} ds \\
\left\| u_{m} [t]-u_{m} [t_{0} ]\right\| _{2}^{2} &\le \left(t-t_{0} \right)\int _{0}^{T}\left\| u_{m,t} [s]\right\| _{2}^{2} ds 
\end{align*}
The above estimates in conjunction with \eqref{GrindEQ__27_} and \eqref{GrindEQ__32_} allow us to conclude that the following estimates hold for all $t,t_{0} \in \left[0,T\right]$: 
\begin{equation} \label{GrindEQ__33_} 
\left\| \varphi _{m} [t]-\varphi _{m} [t_{0} ]\right\| _{2} \le K^{1/2}  \left|t-t_{0} \right|^{1/2}  
\end{equation} 
\begin{equation} \label{GrindEQ__34_} 
\left\| u_{m} [t]-u_{m} [t_{0} ]\right\| _{2} \le \bar{K}^{1/2}  \left|t-t_{0} \right|^{1/2}  
\end{equation}
Thus, by virtue of \eqref{GrindEQ__27_} the sequence $$\left\{\, \varphi _{m} \, :\, m=1,2,...\right\}$$ is bounded in $L^{\infty } \left((0,T); H^{2} (0,1)\right)$ and in $L^{2} \left((0,T);H^{3} (0,1)\right)$. 
Moreover, the sequence $$\left\{\, \varphi _{m,t} \, :\, m=1,2,...\right\}$$ is bounded in $L^{\infty } \left((0,T);L^{2} (0,1)\right)$ and in $L^{2} \left((0,T);H^{1} (0,1)\right)$ while by virtue of \eqref{GrindEQ__32_} the sequence $$\left\{\, u_{m} \, :\, m=1,2,...\right\}$$ is bounded in $L^{\infty } \left((0,T);H_{0}^{1} (0,1)\right)$ and in $L^{2} \left((0,T);H^{2} (0,1)\right)$ and the sequence $$\left\{\, u_{m,t} \, :\, m=1,2,...\right\}$$ is bounded in $L^{2} \left((0,T);L^{2} (0,1)\right)$. Finally, by virtue of \eqref{GrindEQ__27_}, \eqref{GrindEQ__32_} and \eqref{GrindEQ__33_}, \eqref{GrindEQ__34_} the sequences $\left\{\, \varphi _{m} \, :\, m=1,2,...\right\}$, $\left\{\, u_{m} \, :\, m=1,2,...\right\}$ are uniformly bounded in $H^{1} (0,1)$ and are uniformly equicontinuous in $L^{2} (0,1)$. 

Consequently, there exists subsequences 
\begin{align*}
\left\{\, \varphi _{m_{k} } \, :\, k=1,2,...\right\}&\subseteq \left\{\, \varphi _{m} \, :\, m=1,2,...\right\},  \\
\left\{\, u_{m_{k} } \, :\, k=1,2,...\right\} &\subseteq \left\{\, u_{m} \, :\, m=1,2,...\right\} \\
\end{align*}
and
$$\varphi \in L^{\infty } \left((0,T);H^{2} (0,1)\right)\cap L^{2} \left((0,T);H^{3} (0,1)\right),$$ $$ u\in L^{\infty } \left((0,T);H_{0}^{1} (0,1)\right)\cap L^{2} \left((0,T);H^{2} (0,1)\right) $$ with $$ \varphi _{t} \in L^{\infty } \left((0,T);L^{2} (0,1)\right)\cap L^{2} \left((0,T);H^{1} (0,1)\right), u_{t} \in L^{2} \left((0,T);L^{2} (0,1)\right)$$
such that
\begin{align*}
\varphi _{m_{k} } \to \varphi  \textrm{ in } L^{2} \left((0,T);H^{3} (0,1)\right) \textrm{weakly} \\
u_{m_{k} } \to u \textrm{ in } L^{2} \left((0,T);H^{2} (0,1)\right) \textrm{weakly} \\
\varphi _{m_{k} ,t} \to \varphi _{t}  \textrm{ in } L^{2} \left((0,T);H^{1} (0,1)\right) \textrm{weakly} \\
u_{m_{k} ,t} \to u_{t}  \textrm{ in } L^{2} \left((0,T);L^{2} (0,1)\right) \textrm{weakly} \\
\varphi _{m_{k} } \to \varphi  \textrm{ and } u_{m_{k} } \to u \textrm{ in } C^{0} \left([0,T];L^{2} (0,1)\right) \textrm{strongly}
\end{align*}
We show next that 
$$\left(\xi (t),w(t)\right)=\left(\xi _{0} +w_{0} t-\int _{0}^{t}\int _{0}^{\tau }f(s)ds d\tau  ,w_{0} -\int _{0}^{t}f(s)ds \right)$$
for $t\in \left[0,T\right]$ and $\varphi $ constitute a weak solution on $[0,T]$ of the initial-boundary value problem \eqref{GrindEQ__1_}-\eqref{GrindEQ__5_} with initial condition \eqref{GrindEQ__8_}. The fact that
$$u\in C^{0} \left([0,T];H^{1} (0,1)\right) \textrm{ and } \varphi \in C^{0} \left([0,T];H^{2} (0,1)\right)$$ 
follows from Theorem 4 on pages 287-288 in \cite{8}.
Since 
$$\int _{0}^{1}\varphi _{m} (t,x)dx =\int _{0}^{1}\varphi _{m,t} (t,x)dx =0, $$ $$ \varphi _{m,x} (t,1)=\int _{0}^{1}\left(p(x)\varphi _{m,xx} (t,x)+p'(x)\varphi _{m,x} (t,x)\right)dx =0, $$ $$ \varphi _{m,x} (t,0)=\int _{0}^{1}\left(q(x)\varphi _{m,xx} (t,x)+q'(x)\varphi _{m,x} (t,x)\right)dx =0$$
for all $t\in \left[0,T\right]$ and for every $p,q\in C^{1} \left(\left[0,1\right]\right)$ with $p\eqref{GrindEQ__1_}=1$, $p(0)=q\eqref{GrindEQ__1_}=0$, $q(0)=-1$ (recall \eqref{GrindEQ__14_}, \eqref{GrindEQ__15_}), we obtain that 
$$\int _{0}^{T}\int _{0}^{1}v(t)\varphi _{m_{k} } (t,x)dx dt =\int _{0}^{T}\int _{0}^{1}v(t)\varphi _{m_{k} ,t} (t,x)dx dt =0,$$ $$\int _{0}^{T}\int _{0}^{1}v(t)\left(p(x)\varphi _{m_{k} ,xx} (t,x)+p'(x)\varphi _{m_{k} ,x} (t,x)\right)dx dt =0,$$ $$\int _{0}^{T}\int _{0}^{1}v(t)\left(q(x)\varphi _{m_{k} ,xx} (t,x)+q'(x)\varphi _{m_{k} ,x} (t,x)\right)dx dt =0$$ for all $v\in L^{2} (0,T)$. Since $\varphi _{m_{k} } \to \varphi $ in $L^{2} \left((0,T);H^{3} (0,1)\right)$ weakly and $\varphi _{m_{k} ,t} \to \varphi _{t} $ in $L^{2} \left((0,T);H^{1} (0,1)\right)$ weakly we get that 
$$\int _{0}^{T}v(t)\int _{0}^{1}\varphi (t,x)dx dt =\int _{0}^{T}v(t)\int _{0}^{1}\varphi _{t} (t,x)dx dt =0,$$ $$\int _{0}^{T}v(t)\int _{0}^{1}\left(p(x)\varphi _{xx} (t,x)+p'(x)\varphi _{x} (t,x)\right)dx dt =0,$$ $$\int _{0}^{T}v(t)\int _{0}^{1}\left(q(x)\varphi _{xx} (t,x)+q'(x)\varphi _{x} (t,x)\right)dx dt =0$$ for all $v\in L^{2} (0,T)$. Thus, $$\int _{0}^{1}\varphi (t,x)dx =\int _{0}^{1}\varphi _{t} (t,x)dx =0,$$ $$\varphi _{x} (t,1)=\int _{0}^{1}\left(p(x)\varphi _{xx} (t,x)+p'(x)\varphi _{x} (t,x)\right)dx =0,$$ $$\varphi _{x} (t,0)=\int _{0}^{1}\left(q(x)\varphi _{xx} (t,x)+q'(x)\varphi _{x} (t,x)\right)dx =0$$ for a.e. $t\in \left(0,T\right)$. Consequently, definitions \eqref{GrindEQ__6_} and \eqref{GrindEQ__7_} imply that $\varphi \in L^{\infty } \left((0,T);\bar{S}\right)$ and $\varphi _{t} \in L^{\infty } \left((0,T);S\right)$. Since $\varphi \in C^{0} \left([0,T];H^{2} (0,1)\right)$ (which implies that the mappings $t\mapsto \int _{0}^{1}\varphi (t,x)dx $, $t\mapsto \varphi _{x} (t,1)= \int _{0}^{1}\left(p(x)\varphi _{xx} (t,x)\right.$ $\left. +p'(x) \varphi _{x} (t,x)\right)dx $, $t\mapsto \varphi _{x} (t,0)=\int _{0}^{1}\left(q(x)\varphi _{xx} (t,x)+q'(x)\varphi _{x} (t,x)\right)dx $ are continuous), we get that $\varphi \in C^{0} \left([0,T];\bar{S}\right)$. Similarly, we show that $u\in C^{0} \left([0,T];H_{0}^{1} (0,1)\right)$.

Since \eqref{GrindEQ__14_}, \eqref{GrindEQ__15_} and \eqref{GrindEQ__17_} imply that ${\mathop{\lim }\limits_{m\to +\infty }} \left(\left\| \varphi _{0} -\varphi _{m} [0]\right\| _{2} \right)=0$ and since $\varphi _{m_{k} } \to \varphi $ in $C^{0} \left([0,T];L^{2} (0,1)\right)$ strongly we obtain that $\varphi [0]=\varphi _{0} $. 

Since $u_{m} (t,x)=\int _{0}^{x}\varphi _{m,t} (t,s)ds $ for all $t\in \left[0,T\right]$, $x\in \left[0,1\right]$ (recall \eqref{GrindEQ__30_}), we obtain that $\int _{0}^{T}\int _{0}^{1}\left(u_{m} (t,x)-\int _{0}^{x}\varphi _{m,t} (t,s)ds \right)v(t,x)dx dt =0$ for all $v\in L^{2} \left((0,T);L^{2} (0,1)\right)$. Since $\varphi _{m_{k} ,t} \to \varphi _{t} $ in $L^{2} \left((0,T);H^{1} (0,1)\right)$ weakly and $u_{m_{k} } \to u$ in $L^{2} \left((0,T);H^{2} (0,1)\right)$ weakly we get that 
$$\int _{0}^{T}\int _{0}^{1}\left(u(t,x)-\int _{0}^{x}\varphi _{t} (t,s)ds \right)v(t,x)dx dt =0$$ 
for all $v\in L^{2} \left((0,T);L^{2} (0,1)\right)$. Thus $u[t]=\mho \varphi _{t} [t]$ holds for $t\in (0,T)$ a.e. (recall \eqref{GrindEQ__8_}). Consequently, $u_{x} [t]=\varphi _{t} [t]$ for $t\in (0,T)$ a.e.. Since $u\in C^{0} \left([0,T];H_{0}^{1} (0,1)\right)$, it follows that $u_{x} \in C^{0} \left([0,T];L^{2} (0,1)\right)$. Consequently, possibly by redefining $\varphi _{t} $ on a measure zero set, we may assume that $\varphi _{t} \in C^{0} \left([0,T];L^{2} (0,1)\right)$ and that \eqref{GrindEQ__10_} holds. 

Since $u_{m_{k} } \to u$ in $C^{0} \left([0,T];L^{2} (0,1)\right)$ strongly we get $${\mathop{\lim }\limits_{k\to +\infty }} \left(\left\| u_{m_{k} } [0]-u[0]\right\| _{2} \right)=0.$$ Definition \eqref{GrindEQ__28_} in conjunction with \eqref{GrindEQ__17_}, the fact that the functions $g_{n} $ defined by \eqref{GrindEQ__29_} constitute an orthonormal basis of $L^{2} (0,1)$ and Parseval's identity imply that $$\left\| u_{m} [0]-\mho \bar{\varphi }_{0} \right\| _{2}^{2} =\sum _{n=m+1}^{\infty }\left(\mho \bar{\varphi }_{0} ,g_{n} \right)^{2}. $$ Consequently, we get $${\mathop{\lim }\limits_{m\to +\infty }} \left(\left\| u_{m} [0]-\mho \bar{\varphi }_{0} \right\| _{2} \right)=0$$ which gives $u[0]=\mho \bar{\varphi }_{0} $. Using \eqref{GrindEQ__10_} we get that $\varphi _{t} [0]=\bar{\varphi }_{0} $.

Exploiting \eqref{GrindEQ__30_} and \eqref{GrindEQ__31_} we get for every $v\in L^{2} \left((0,T);L^{2} (0,1)\right)$: 
\begin{align} \label{GrindEQ__35_} 
&\int _{0}^{T}\int _{0}^{1}\left(u_{m,t} (t,x)-\varphi _{m,x} (t,x)+\kappa u_{m} (t,x)-\mu u_{m,xx} (t,x)\right)v(t,x)dxdt   \nonumber \\ 
&+\int _{0}^{T}\int _{0}^{1}\left(  \sigma \varphi _{m,xxx} (t,x)+f(t)  \right)           v(t,x)dxdt \nonumber \\
&-\int _{0}^{T}\int _{0}^{1}f(t)\left(1+\sqrt{2} \sum _{n=1}^{m}\frac{(-1)^{n} -1}{n\pi } g_{n} (x) \right)v(t,x)dxdt  =0  
\end{align}
Taking into account the fact that the functions $g_{n} $ defined by \eqref{GrindEQ__29_} constitute an orthonormal basis of $L^{2} (0,1)$ and Parseval's identity, we get 
\[\left\| \chi _{[0,1]} -\sum _{n=1}^{m}\frac{\left(\phi _{n} (0)-\phi _{n} (1)\right)}{n\pi } g_{n}  \right\| _{2}^{2} =\sum _{n=m+1}^{\infty }\frac{\left(\phi _{n} (0)-\phi _{n} (1)\right)^{2} }{n^{2} \pi ^{2} }  \le 4\sum _{n=m+1}^{\infty }\frac{1}{n^{2} \pi ^{2} }  \] 
which combined with \eqref{GrindEQ__15_} and the Cauchy-Schwarz inequality gives:
\begin{align} \label{GrindEQ__36_} 
\left|\int _{0}^{T}\int _{0}^{1}f(t)\left(1+\sqrt{2} \sum _{n=1}^{m}\frac{(-1)^{n} -1}{n\pi } g_{n} (x) \right)v(t,x)dxdt  \right| \nonumber \\ 
\le 2\left(\sum _{n=m+1}^{\infty }\frac{1}{n^{2} \pi ^{2} }  \right)^{1/2} \left(\int _{0}^{T}\left|f(t)\right|^{2} dt \right)^{1/2} \left(\int _{0}^{T}\left\| v[t]\right\| _{2}^{2} dt \right)^{1/2}  
\end{align}
Therefore, it follows from \eqref{GrindEQ__36_} that $${\mathop{\lim }\limits_{m\to +\infty }} \left(\int _{0}^{T}\int _{0}^{1}f(t)\left(1+\sqrt{2} \sum _{n=1}^{m}\frac{(-1)^{n} -1}{n\pi } g_{n} (x) \right)v(t,x)dxdt  \right)=0.$$
Since $\varphi _{m_{k} } \to \varphi $ in $L^{2} \left((0,T);H^{3} (0,1)\right)$ weakly and $u_{m_{k} } \to u$ in $L^{2} \left((0,T); \right.$ $\left. H^{2} (0,1)\right)$ weakly, we get from \eqref{GrindEQ__35_}:
\begin{align} \label{GrindEQ__37_} 
&\int _{0}^{T}\int _{0}^{1}\left(u_{t} (t,x)-\varphi _{x} (t,x)+\kappa u(t,x)-\mu u_{xx} (t,x)\right)v(t,x)dxdt   \nonumber \\
&+\int _{0}^{T}\int _{0}^{1}\left( \sigma \varphi _{xxx} (t,x)+f(t) \right)v(t,x)dxdt=0
\end{align}
Notice that $u_{t} [t]-\varphi _{x} [t]+\kappa u[t]-\mu u_{xx} [t]+\sigma \varphi _{xxx} [t]+f(t)\chi _{[0,1]} $ is of class $L^{2} \left((0,T);L^{2} (0,1)\right)$. Since \eqref{GrindEQ__37_} holds for every $v\in L^{2} \left((0,T);L^{2} (0,1)\right)$, we obtain \eqref{GrindEQ__11_}.

\noindent \underbar{2${}^{nd}$ Part: Proof of \eqref{GrindEQ__13_}.} It follows from the regularity properties of the weak solution that the mapping $\left[0,T\right]\backepsilon t\mapsto W\left(\varphi [t],u[t]\right)$, with $W$ being defined by \eqref{GrindEQ__12_}, is absolutely continuous and satisfies the following equation for a.e. $t\in \left(0,T\right)$: 
\begin{align} \label{GrindEQ__38_} 
 \frac{d}{d\, t} W\left(\varphi [t],u[t]\right)=&\left(u[t],u_{t} [t]\right)+(1+r)\left(\varphi [t],\varphi _{t} [t]\right)+\sigma (1+r)\left(\varphi _{x} [t],\varphi _{tx} [t]\right) \nonumber \\ 
&+r\left(u[t]-\mu \varphi _{x} [t]+\kappa \mho \varphi [t],u_{t} [t]-\mu \varphi _{tx} [t]+\kappa \left(\mho \varphi [t]\right)_{t} \right) 
\end{align} 
Using \eqref{GrindEQ__8_}, \eqref{GrindEQ__10_}, \eqref{GrindEQ__11_} (which also imply the equalities $\left(\mho \varphi [t]\right)_{t} =u[t]$, $\varphi _{tx} [t]=u_{xx} [t]$, $\left(\mho \varphi [t]\right)_{x} =\varphi [t]$) we obtain from \eqref{GrindEQ__38_} for a.e. $t\in \left(0,T\right)$:
\begin{align} \label{GrindEQ__39_} 
\frac{d}{d\, t} W\left(\varphi [t],u[t]\right)=&-\kappa \left\| u[t]\right\| _{2}^{2} -\mu \left\| \varphi _{t} [t]\right\| _{2}^{2}  \nonumber \\ 
&-\kappa r\left\| \varphi [t]\right\| _{2}^{2} -r\left(\sigma \kappa +\mu \right)\left\| \varphi _{x} [t]\right\| _{2}^{2} -\mu \sigma r\left\| \varphi _{xx} [t]\right\| _{2}^{2} \nonumber \\ 
&-f(t)\left(r\kappa \int _{0}^{1}\left(\mho \varphi [t]\right)(x)dx \right) \nonumber \\
&-f(t) \left( (r+1)\int _{0}^{1}u(t,x)dx -r\mu \left(\varphi (t,1)-\varphi (t,0)\right) \right)
\end{align}
Integrating equation \eqref{GrindEQ__39_} and using \eqref{GrindEQ__8_} (which implies that $\int _{0}^{1}\left(\mho \varphi [t]\right)(x)dx =-\int _{0}^{1}x\varphi (t,x)dx $) and \eqref{GrindEQ__10_} we obtain \eqref{GrindEQ__13_}. 

\noindent \underbar{3${}^{rd}$ Part: Uniqueness.} We next exploit \eqref{GrindEQ__13_} in order to prove that the weak solution on $[0,T]$ of the initial-boundary value problem \eqref{GrindEQ__1_}-\eqref{GrindEQ__5_} with initial condition \eqref{GrindEQ__9_} is unique.

The difference of two weak solutions on $[0,T]$ of the initial-boundary value problem \eqref{GrindEQ__1_}-\eqref{GrindEQ__5_} with initial condition \eqref{GrindEQ__9_} is a weak solution on $[0,T]$ of the initial-boundary value problem \eqref{GrindEQ__1_}-\eqref{GrindEQ__5_} with $\left(\xi _{0} ,w_{0} \right)=\left(0,0\right)$, $f(t)\equiv 0$, $\varphi _{0} =\bar{\varphi }_{0} =0$. We next show that the weak solution on $[0,T]$ of the initial-boundary value problem \eqref{GrindEQ__1_}-\eqref{GrindEQ__5_} with $\left(\xi _{0} ,w_{0} \right)=\left(0,0\right)$, $f(t)\equiv 0$, $\varphi _{0} =\bar{\varphi }_{0} =0$ is necessarily identically zero.

Using \eqref{GrindEQ__13_} and definition \eqref{GrindEQ__12_} with $r=1$ we get that $W\left(\varphi [t],u[t]\right)=0$ for all $t\in [0,T]$. Definition \eqref{GrindEQ__12_} implies that $\varphi [t]=0$ for all $t\in [0,T]$. Moreover, the unique solution of \eqref{GrindEQ__1_} with $\left(\xi _{0} ,w_{0} \right)=\left(0,0\right)$, $f(t)\equiv 0$ is the solution $\left(\xi (t),w(t)\right)=\left(0,0\right)$ (uniqueness of solutions for linear ODEs). The proof is complete. 
\qed
\end{proof}  

The proof of Lemma \ref{lemma1} is given next. 
 
\begin{proof}{\hspace{-0.15em}\textit{of Lemma \ref{lemma1}}} 
Since
\begin{align} \label{GrindEQ__X1_}
\inf & \left\lbrace  \frac{1}{1-\bar{G}p}+\frac{\left(r+1\right)^{2}\left(1+\sigma \pi ^{2}\right) \left(1+p \right)}{p \mu r \left(\kappa + \mu \pi ^{2} \right)} 
 : p>0, \bar{G}p <1  \right\rbrace \nonumber \\
&=\left(1+\frac{\kappa}{\pi \mu \sqrt{3}} \right)^{2}+\frac{\left( r+1\right)^{2}\left(1+\sigma \pi ^{2} \right)}{r \mu \left(\kappa + \mu \pi ^{2} \right)}
\end{align}
where
\begin{align*}
\bar{G}:=\frac{\kappa ^{2} r \left(\kappa + \mu \pi ^{2} \right)}{3 \mu \pi^{2} \left(r+1 \right)^{2}\left(1+\sigma \pi ^{2}\right)}
\end{align*}
it follows from the fact that $\gamma>\gamma^{*}$ and definition \eqref{GrindEQ__A1_} of $\gamma^{*}$ that there exists $p>0$ with $\bar{G}p <1$ such that
\begin{align} \label{GrindEQ__X2_}
&\gamma\left(1-\frac{\mu r k^{3}Q}{1+\sigma \pi ^{2}}\left( \frac{1}{\displaystyle{1-\frac{\kappa ^{2}rp \left(\kappa + \mu \pi ^{2}\right)}{3 \mu \pi^{2} \left(r+1\right)^{2}\left(1+\sigma \pi ^{2}\right)}}}+\frac{\left(r+1\right)^{2}\left(1+\sigma \pi ^{2}\right) }{\mu r \left(\kappa + \mu \pi ^{2} \right)} \frac{\left(1+p \right)}{p } \right) \right) \nonumber \\
&>2
\end{align}
Define
\begin{equation} \label{GrindEQ__43_} 
\delta :=\frac{\left(1+\sigma \pi ^{2} \right)(r+1)^{2} \left(1+p\right)}{p\left(\kappa +\mu \pi ^{2} \right)r\mu } -\frac{\kappa ^{2} \left(1+p\right)}{3\mu ^{2} \pi ^{2} }  
\end{equation} 
The fact that $3\mu \pi ^{2} \left(1+\sigma \pi ^{2} \right)(r+1)^{2} >\kappa ^{2} p\left(\kappa +\mu \pi ^{2} \right)r$ guarantees that $\delta >0$. Definition \eqref{GrindEQ__43_} and inequality \eqref{GrindEQ__X2_} guarantee that the following inequalities are valid:
\begin{align} \label{GrindEQ__44_} 
\gamma &>\frac{\gamma ^{2} k^{3} r\left(1+\delta \right)Q\left(3\delta \mu ^{2} \pi ^{2} +\kappa ^{2} \left(1+p\right)\right)}{6\delta \mu \pi ^{2} \left(1+\sigma \pi ^{2} \right)+\gamma k^{3} r\left(1+\delta \right)Q\left(3\delta \mu ^{2} \pi ^{2} +\kappa ^{2} \left(1+p\right)\right)} +2 \nonumber  \\ 
\gamma &>\frac{\gamma ^{2} k^{3} (r+1)^{2} Q\left(1+\delta \right)\left(1+p\right)}{2\delta p\left(\kappa +\mu \pi ^{2} \right)+\gamma k^{3} (r+1)^{2} Q\left(1+\delta \right)\left(1+p\right)} +2
\end{align} 
Inequalities \eqref{GrindEQ__44_} guarantee the existence of $\varepsilon \in \left(0,1\right)$ such that 
\begin{align} \label{GrindEQ__45_} 
 1-\frac{2 }{\gamma } >\varepsilon >&\frac{\gamma k^{3} r\left(1+\delta \right)Q\left(3\delta \mu ^{2} \pi ^{2} +\kappa ^{2} \left(1+p\right)\right)}{6\delta \mu \pi ^{2} \left(1+\sigma \pi ^{2} \right)+\gamma k^{3} r\left(1+\delta \right)Q\left(3\delta \mu ^{2} \pi ^{2} +\kappa ^{2} \left(1+p\right)\right)}, \nonumber \\ 
 1-\frac{2 }{\gamma } >\varepsilon >&\frac{\gamma k^{3} (r+1)^{2} Q\left(1+\delta \right)\left(1+p\right)}{2\delta p\left(\kappa +\mu \pi ^{2} \right)+\gamma k^{3} (r+1)^{2} Q\left(1+\delta \right)\left(1+p\right)}  
\end{align} 
Inequalities \eqref{GrindEQ__45_} allow us to derive the following inequalities:
\begin{align} \label{GrindEQ__46_} 
 1+\sigma \pi ^{2} &>\frac{\gamma k^{3} r(1-\varepsilon )}{2\varepsilon } \left(1+\delta \right)Q\left(\mu +\frac{\kappa ^{2} \left(1+p\right)}{3\delta \mu \pi ^{2} } \right) \nonumber \\ 
 \kappa +\mu \pi ^{2} &>\frac{\gamma k^{3} (r+1)^{2} (1-\varepsilon )}{2\varepsilon \delta p} Q\left(1+\delta \right)\left(1+p\right) \nonumber \\ 
 \gamma (1-\varepsilon )&>2  
\end{align} 
Clearly, we conclude from \eqref{GrindEQ__46_} that there exists $\lambda \in \left(0,1\right)$ (sufficiently close to 1) such that the following inequalities hold: 
\begin{align} \label{GrindEQ__47_} 
 1+\lambda \sigma \pi ^{2} &>\frac{\gamma k^{3} r(1-\varepsilon )}{2\varepsilon } \left(1+\delta \right)Q\left(\mu +\frac{\kappa ^{2} \left(1+p\right)}{3\delta \mu \pi ^{2} } \right) \nonumber \\ 
 \kappa +\lambda \mu \pi ^{2} &>\frac{\gamma k^{3} (r+1)^{2} (1-\varepsilon )}{2\varepsilon \delta p} Q\left(1+\delta \right)\left(1+p\right) \nonumber \\ 
 \gamma  (1-\varepsilon )&> 2  
\end{align}
Let arbitrary $\varphi \in \bar{S}$, $u\in H_{0}^{1} (0,1)$, $\left(\xi ,w\right)\in {\mathbb R}^{2} $ be given. Then the following inequality holds: 
\begin{align} \label{GrindEQ__48_} 
&k Q r  \left(w+k\xi \right)\int _{0}^{1}\left(\kappa x\varphi (x)-\frac{r+1}{r} u(x)+\mu \varphi '(x)\right)dx  \\ 
&\le  \frac{\varepsilon }{2k } \left(w+k\xi \right)^{2} +\frac{k^{3} r^{2} }{2\varepsilon } Q^{2} \left(\int _{0}^{1}\left(\kappa x\varphi (x)-\frac{r+1}{r} u(x)+\mu \varphi '(x)\right)dx \right)^{2}  
\end{align}
Combining \eqref{GrindEQ__40_} and \eqref{GrindEQ__48_} we get:
\begin{align} \label{GrindEQ__49_} 
&q\left(\xi ,w,\varphi ,u\right) \nonumber \\
\ge& k\xi ^{2} +\frac{\gamma  (1-\varepsilon )-2 }{2k } \left(w+k\xi \right)^{2} +Q\mu \left\| u'\right\| _{2}^{2} \nonumber \\
&+Q\mu r\left\| \varphi '\right\| _{2}^{2} +Q\mu \sigma r\left\| \varphi ''\right\| _{2}^{2} +Q\kappa \left\| u\right\| _{2}^{2} \nonumber \\ 
&-\frac{\gamma k^{3} r^{2} (1-\varepsilon )}{2\varepsilon } Q^{2} \left(\int _{0}^{1}\left(\kappa x\varphi (x)-\frac{r+1}{r} u(x)\right)dx+\mu \left(\varphi (1)-\varphi (0)\right) \right)^{2} 
\end{align} 
Moreover, the following inequality holds: 
\begin{align} \label{GrindEQ__50_} 
&\mu \left(\varphi (1)-\varphi (0)\right)\int _{0}^{1}\left(\kappa x\varphi (x)-\frac{r+1}{r} u(x)\right)dx \nonumber  \\ 
&\le \frac{\delta }{2} \mu ^{2} \left(\varphi (1)-\varphi (0)\right)^{2} +\frac{1}{2\delta } \left(\int _{0}^{1}\left(\kappa x\varphi (x)-\frac{r+1}{r} u(x)\right)dx \right)^{2} 
\end{align}
Combining \eqref{GrindEQ__49_} and \eqref{GrindEQ__50_} we get:
\begin{align} \label{GrindEQ__51_} 
q\left(\xi ,w,\varphi ,u\right)\ge & k\xi ^{2} +\frac{\gamma (1-\varepsilon )-2 }{2k } \left(w+k\xi \right)^{2} +Q\mu \left\| u'\right\| _{2}^{2}  \nonumber \\ 
&+Q\mu r\left\| \varphi '\right\| _{2}^{2} +Q\mu \sigma r\left\| \varphi ''\right\| _{2}^{2} +Q\kappa \left\| u\right\| _{2}^{2}  \nonumber \\ & -\frac{\gamma k^{3} r^{2} (1-\varepsilon )}{2\varepsilon \delta } Q^{2} \left(1+\delta \right)\left(\int _{0}^{1}\left(\kappa x\varphi (x)-\frac{r+1}{r} u(x)\right)dx \right)^{2}  \nonumber \\ 
&-\frac{\gamma k^{3} r^{2} (1-\varepsilon )}{2\varepsilon } \left(1+\delta \right)Q^{2} \mu ^{2} \left(\varphi (1)-\varphi (0)\right)^{2}  
\end{align}
Using the Cauchy-Schwarz inequality we obtain the estimate $\left|\varphi \eqref{GrindEQ__1_}-\varphi (0)\right|\le \left\| \varphi '\right\| _{2} $. Therefore, we get from \eqref{GrindEQ__51_}:
\begin{align}  \label{GrindEQ__52_}
q\left(\xi ,w,\varphi ,u\right)\ge& k\xi ^{2} +\frac{\gamma (1-\varepsilon )-2 }{2k} \left(w+k\xi \right)^{2}  \nonumber \\  
&+Q\mu \left\| u'\right\| _{2}^{2} +Q\kappa \left\| u\right\| _{2}^{2}  \nonumber \\
&+Q\mu r\left(1-\frac{\gamma k^{3} r(1-\varepsilon )}{2\varepsilon } \left(1+\delta \right)Q\mu \right)\left\| \varphi '\right\| _{2}^{2} +Q\mu \sigma r\left\| \varphi ''\right\| _{2}^{2} \nonumber \\ 
&-\frac{\gamma k^{3} r^{2} (1-\varepsilon )}{2\varepsilon \delta } Q^{2} \left(1+\delta \right)\left(\kappa \int _{0}^{1}x\varphi (x)dx -\frac{r+1}{r} \int _{0}^{1}u(x)dx \right)^{2} 
\end{align}  
Using the inequality $\displaystyle{\left(a+b\right)^{2} \le \left(1+p\right)a^{2} +\frac{1+p}{p} b^{2} }$ that holds for all $a,b\in {\mathbb R}$, we get from \eqref{GrindEQ__52_}:
\begin{align} \label{GrindEQ__53_} 
q\left(\xi ,w,\varphi ,u\right)\ge & k\xi ^{2} +\frac{\gamma (1-\varepsilon )-2 }{2k } \left(w+k\xi \right)^{2} +Q\mu \left\| u'\right\| _{2}^{2} \nonumber \\
&+Q\kappa \left\| u\right\| _{2}^{2} +Q\mu r\left(1-\frac{\gamma k^{3} r(1-\varepsilon )}{2\varepsilon } \left(1+\delta \right)Q\mu \right)\left\| \varphi '\right\| _{2}^{2}  \nonumber \\
&+Q\mu \sigma r\left\| \varphi ''\right\| _{2}^{2}\nonumber \\ 
&-\frac{\gamma k^{3} r^{2} (1-\varepsilon )}{2\varepsilon \delta } Q^{2} \kappa ^{2} \left(1+\delta \right)\left(1+p\right)\left(\int _{0}^{1}x\varphi (x)dx \right)^{2}  \nonumber \\
&-\frac{\gamma k^{3} (r+1)^{2} (1-\varepsilon )}{2\varepsilon \delta p} Q^{2} \left(1+\delta \right)\left(1+p\right)\left(\int _{0}^{1}u(x)dx \right)^{2}   
\end{align}
Using the Cauchy-Schwarz inequality we get the inequalities $\left|\int _{0}^{1}x\varphi (x)dx \right|\le \frac{1}{\sqrt{3} } \left\| \varphi \right\| _{2} $ and $\left|\int _{0}^{1}u(x)dx \right|\le \left\| u\right\| _{2} $. Therefore, we obtain from \eqref{GrindEQ__53_}:
\begin{align} \label{GrindEQ__54_} 
q\left(\xi ,w,\varphi ,u\right)\ge& k\xi ^{2} +\frac{\gamma (1-\varepsilon )-2 }{2k } \left(w+k\xi \right)^{2} +Q\mu \left\| u'\right\| _{2}^{2} \nonumber \\ 
&+Q\mu r\left(1-\frac{\gamma k^{3} r(1-\varepsilon )}{2\varepsilon } \left(1+\delta \right)Q\mu \right)\left\| \varphi '\right\| _{2}^{2} +Q\mu \sigma r\left\| \varphi ''\right\| _{2}^{2} \nonumber \\ 
&-\frac{\gamma k^{3} r^{2} (1-\varepsilon )}{6\varepsilon \delta } Q^{2} \kappa ^{2} \left(1+\delta \right)\left(1+p\right)\left\| \varphi \right\| _{2}^{2}  \nonumber \\
&+Q\kappa \left\| u\right\| _{2}^{2} -\frac{\gamma k^{3} (r+1)^{2} (1-\varepsilon )}{2\varepsilon \delta p} Q^{2} \left(1+\delta \right)\left(1+p\right)\left\| u\right\| _{2}^{2}   
\end{align} 
Since $u,\varphi '\in H_{0}^{1} (0,1)$ and $\int _{0}^{1}\varphi (x)dx =0$ (recall definition \eqref{GrindEQ__6_} and the fact that $\varphi \in \bar{S}$), we obtain from Wirtinger's inequality the estimates $\left\| u'\right\| _{2}^{2} \ge \pi ^{2} \left\| u\right\| _{2}^{2} $, $\left\| \varphi '\right\| _{2}^{2} \ge \pi ^{2} \left\| \varphi \right\| _{2}^{2} $, $\left\| \varphi ''\right\| _{2}^{2} \ge \pi ^{2} \left\| \varphi '\right\| _{2}^{2} $. Consequently, we get from \eqref{GrindEQ__54_}:
\begin{align}  \label{GrindEQ__55_}
&q\left(\xi ,w,\varphi ,u\right) \nonumber \\
\ge& k\xi ^{2} +\frac{\gamma (1-\varepsilon )-2 }{2k } \left(w+k\xi \right)^{2} +\left(1-\lambda \right)Q\mu \left\| u'\right\| _{2}^{2} \nonumber \\ 
&+Q\mu r\left(1+\lambda \sigma \pi ^{2} -\frac{\gamma k^{3} r(1-\varepsilon )}{2\varepsilon } \left(1+\delta \right)Q\left(\mu +\frac{\kappa ^{2} \left(1+p\right)}{3\delta \mu \pi ^{2} } \right)\right)\left\| \varphi '\right\| _{2}^{2} \nonumber \\ 
&+\left(1-\lambda \right)Q\mu \sigma r\left\| \varphi ''\right\| _{2}^{2}  \nonumber \\
&+Q\left(\kappa +\lambda \mu \pi ^{2} -\frac{\gamma k^{3} (r+1)^{2} (1-\varepsilon )}{2\varepsilon \delta p} Q\left(1+\delta \right)\left(1+p\right)\right)\left\| u\right\| _{2}^{2}
\end{align} 
The existence of a constant $A>0$ for which \eqref{GrindEQ__42_} holds is a consequence of inequality \eqref{GrindEQ__55_} and inequalities \eqref{GrindEQ__47_} as well as the inequality $\left\| \varphi '\right\| _{2}^{2} \ge \pi ^{2} \left\| \varphi \right\| _{2}^{2} $. The proof is complete. 
\qed
\end{proof}
The proof of Theorem \ref{theorem2} is provided next. 
\begin{proof}{\hspace{-0.15em}\textit{of Theorem \ref{theorem2}}}
Define $\varphi _{m} $ and $u_{m} $ by means of \eqref{GrindEQ__14_} and \eqref{GrindEQ__28_}, respectively, where $\phi _{n} ,g_{n} $ are given by \eqref{GrindEQ__15_}, \eqref{GrindEQ__29_} and all $a_{n} $ for $n=1,...,m$ are given by the solution of the following initial value problem \eqref{GrindEQ__17_} with
\begin{equation} \label{GrindEQ__62_} 
\ddot{a}_{n} (t)=-\left(\mu n^{2} \pi ^{2} +\kappa \right)\dot{a}_{n} (t)-\left(\sigma n^{2} \pi ^{2} +1\right)n^{2} \pi ^{2} a_{n} (t)+\beta _{n} f_{m} (t) 
\end{equation} 
\begin{equation} \label{GrindEQ__63_} 
\dot{\xi }_{m} (t)=w_{m} (t)\quad ,\quad \dot{w}_{m} (t)=-f_{m} (t) 
\end{equation} 
\begin{equation} \label{GrindEQ__64_} 
\left(\xi _{m} (0),w_{m} (0)\right)=\left(\xi _{0} ,w_{0} \right) 
\end{equation} 
where 
\begin{equation} \label{GrindEQ__65_} 
\beta _{n} =\sqrt{2} \left((-1)^{n} -1\right) 
\end{equation} 
\begin{align} \label{GrindEQ__66_} 
f_{m} (t)=&\gamma k \left(w_{m} (t)+k\xi _{m} (t)\right) \nonumber \\
&-\gamma k^{3} rQ\sum _{n=1}^{m}\beta _{n} \frac{\left(\kappa +\mu n^{2} \pi ^{2} \right)a_{n} (t)}{n^{2} \pi ^{2} }  -\gamma k^{3} (r+1)Q\sum _{n=1}^{m}\beta _{n} \frac{\dot{a}_{n} (t)}{n^{2} \pi ^{2} }   \nonumber \\ 
=&\gamma k \left(w_{m} (t)+k\xi _{m} (t)\right)
\nonumber \\
&-\gamma k^{3} r Q\int _{0}^{1}\left(\kappa x\varphi _{m} (t,x)-\frac{r+1}{r} u_{m} (t,x)+\mu \varphi _{m,x} (t,x)\right)dx 
\end{align} 
\noindent Furthermore, define for $(t,x)\in \left[0,T\right]\times \left[0,1\right]$:
\begin{equation} \label{GrindEQ__67_} 
\theta _{m} (t,x):=\sum _{n=1}^{m}\frac{a_{n} (t)}{n\pi } g_{n} (x) =\int _{0}^{x}\varphi _{m} (t,s)ds  
\end{equation} 
\begin{equation} \label{GrindEQ__68_} 
V_{m} (t):=\frac{1}{2} \xi _{m}^{2} (t)+\frac{1}{2k^{2} } \left(w_{m} (t)+k\xi _{m} (t)\right)^{2} +QW_{m} (t) 
\end{equation} 
\begin{align} \label{GrindEQ__69_}
W_{m} (t):=&\frac{1}{2} \left\| u_{m} [t]\right\| _{2}^{2} +\frac{1+r}{2} \left\| \varphi _{m} [t]\right\| _{2}^{2} \nonumber \\
&+\frac{\sigma (1+r)}{2} \left\| \varphi _{m,x} [t]\right\| _{2}^{2} +\frac{r}{2} \left\| u_{m} [t]-\mu \varphi _{m,x} [t]+\kappa \theta _{m} [t]\right\| _{2}^{2} \nonumber  \\  
=&\frac{1}{2} \sum _{n=1}^{m}\frac{\dot{a}_{n}^{2} (t)}{n^{2} \pi ^{2} }  +\frac{1+r}{2} \sum _{n=1}^{m}a_{n}^{2} (t) +\frac{\sigma (1+r)}{2} \sum _{n=1}^{m}n^{2} \pi ^{2} a_{n}^{2} (t) \nonumber \\
&+\frac{r}{2} \sum _{n=1}^{m}\frac{\left(\dot{a}_{n} (t)+\left(\mu n^{2} \pi ^{2} +\kappa \right)a_{n} (t)\right)^{2} }{n^{2} \pi ^{2} }   
\end{align} 
Using \eqref{GrindEQ__68_}, \eqref{GrindEQ__69_} and \eqref{GrindEQ__62_}, \eqref{GrindEQ__63_}, \eqref{GrindEQ__66_} as well as definition \eqref{GrindEQ__40_} we obtain for all $t\in \left[0,T\right]$:
\begin{equation} \label{GrindEQ__70_} 
\dot{V}_{m} (t)=-q\left(\xi _{m} (t),w_{m} (t),\varphi _{m} [t],u_{m} [t]\right)-\frac{1}{2\gamma k^{3} } f_{m}^{2} (t) 
\end{equation}
Definitions \eqref{GrindEQ__68_}, \eqref{GrindEQ__69_} as well as the fact $\pi ^{2} \left\| \theta _{m} [t]\right\| _{2}^{2} \le \left\| \varphi _{m} [t]\right\| _{2}^{2} $ (Wirtinger's inequality; recall definition \eqref{GrindEQ__67_}) imply the existence of a constant $c_{1} >0$ (independent of the particular solution and $m$, $T>0$) for which the following estimate holds for all $t\in \left[0,T\right]$:
\begin{equation} \label{GrindEQ__71_} 
V_{m} (t)\le c_{1} \left(\xi _{m}^{2} (t)+w_{m}^{2} (t)+\left\| u_{m} [t]\right\| _{2}^{2} +\left\| \varphi _{m} [t]\right\| _{H^{1} (0,1)}^{2} \right) 
\end{equation} 
Exploiting \eqref{GrindEQ__42_} and \eqref{GrindEQ__71_} we conclude that the differential inequality $\dot{V}_{m} (t)\le -c_{1}^{-1} AV_{m} (t)$ holds for all $t\in \left[0,T\right]$. Consequently, we get:
\begin{equation} \label{GrindEQ__72_}
V_{m} (t)\le \exp \left(-c_{1}^{-1} At\right)V_{m} (0), \textrm{for all } t\in \left[0,T\right]
\end{equation}                                       
Definitions \eqref{GrindEQ__66_}, \eqref{GrindEQ__68_}, \eqref{GrindEQ__69_} imply the existence of a constant $c_{2} >0$ (independent of the particular solution and $m$, $T>0$) for which the following estimate holds for all $t\in \left[0,T\right]$:
\begin{equation} \label{GrindEQ__73_} 
f_{m}^{2} (t)\le c_{2} V_{m} (t) 
\end{equation}
Moreover, equations \eqref{GrindEQ__14_}, \eqref{GrindEQ__17_} and \eqref{GrindEQ__28_} imply that $\left\| u_{m} [0]\right\| _{2} \le \left\| \varphi _{m,t} [0]\right\| _{2} $, $\left\| \varphi _{m,x} [0]\right\| _{2} \le \left\| \varphi '_{0} \right\| _{2} $, $\left\| \varphi _{m} [0]\right\| _{2} \le \left\| \varphi _{0} \right\| _{2} $ and $\left\| \varphi _{m,t} [0]\right\| _{2} \le \left\| \bar{\varphi }_{0} \right\| _{2} $. Consequently, we obtain from \eqref{GrindEQ__71_} and \eqref{GrindEQ__64_}: 
\begin{equation} \label{GrindEQ__74_} 
V_{m} (0)\le c_{1} \left(\xi _{0}^{2} +w_{0}^{2} +\left\| \bar{\varphi }_{0} \right\| _{2}^{2} +\left\| \varphi _{0} \right\| _{H^{1} (0,1)}^{2} \right) 
\end{equation} 
Working exactly as in the proof of Theorem \ref{theorem1} we obtain the following estimate (the analogue of \eqref{GrindEQ__26_}) for $t\in \left[0,T\right]$:
\begin{align} \label{GrindEQ__75_}  
&\left\| \varphi _{m,t} [t]-\mu \varphi _{m,xx} [t]+\kappa \varphi _{m} [t]\right\| _{2}^{2} +2\left\| \varphi _{m,t} [t]\right\| _{2}^{2}   \nonumber \\ 
&+3\left\| \varphi _{m,x} [t]\right\| _{2}^{2}+3\sigma \left\| \varphi _{m,xx} [t]\right\| _{2}^{2}\nonumber \\ 
&+\int _{0}^{t}\left(\mu \sigma \left\| \varphi _{m,xxx} [s]\right\| _{2}^{2} +\left(\kappa \sigma +\mu \right)\left\| \varphi _{m,xx} [s]\right\| _{2}^{2} +\kappa \left\| \varphi _{m,x} [s]\right\| _{2}^{2} \right)ds  \nonumber \\
&+\int _{0}^{t}\left(\mu \left\| \varphi _{m,tx} [s]\right\| _{2}^{2} +\kappa \left\| \varphi _{m,t} [s]\right\| _{2}^{2} \right)ds  \nonumber \\
\le& \left\| \varphi _{m,t} [0]-\mu \varphi _{m,xx} [0]+\kappa \varphi _{m} [0]\right\| _{2}^{2} \nonumber \\
&+3\left\| \varphi _{m,x} [0]\right\| _{2}^{2} +3\sigma \left\| \varphi _{m,xx} [0]\right\| _{2}^{2} +2\left\| \varphi _{m,t} [0]\right\| _{2}^{2}   \nonumber \\
&+\left(\frac{\left(\mu +\kappa \right)^{2} }{\sigma } +3\right)\left(\sum _{n=1}^{m}\frac{8}{\mu n^{2} \pi ^{2} +\kappa }  \right)\int _{0}^{t}f_{m}^{2} (s)ds
\end{align}
Using \eqref{GrindEQ__72_}, \eqref{GrindEQ__73_}, \eqref{GrindEQ__74_} and the facts that $\left\| \varphi _{m,x} [0]\right\| _{2} \le \left\| \varphi '_{0} \right\| _{2} $, $\left\| \varphi _{m} [0]\right\| _{2} \le \left\| \varphi _{0} \right\| _{2} $, $\left\| \varphi _{m,xx} [0]\right\| _{2} \le \left\| \varphi ''_{0} \right\| _{2} $ and $\left\| \varphi _{m,t} [0]\right\| _{2} \le \left\| \bar{\varphi }_{0} \right\| _{2} $ (recall \eqref{GrindEQ__14_}, \eqref{GrindEQ__17_}) we get from \eqref{GrindEQ__75_} for $t\in \left[0,T\right]$: 
\begin{align} \label{GrindEQ__76_} 
&2\left\| \varphi _{m,t} [t]\right\| _{2}^{2} +3\left\| \varphi _{m,x} [t]\right\| _{2}^{2} +3\sigma \left\| \varphi _{m,xx} [t]\right\| _{2}^{2} \nonumber \\ 
&+\mu \int _{0}^{t}\left(\sigma \left\| \varphi _{m,xxx} [s]\right\| _{2}^{2} +\left\| \varphi _{m,xx} [s]\right\| _{2}^{2} +\left\| \varphi _{m,tx} [s]\right\| _{2}^{2} \right)ds \nonumber \\ \le & 4\left(\sigma +\mu ^{2} +\kappa ^{2} +1\right)\left\| \varphi _{0} \right\| _{H^{2} (0,1)}^{2} +4\left\| \bar{\varphi }_{0} \right\| _{2}^{2} \nonumber \\ 
&+\mu ^{-1} A^{-1} c_{1}^{2} c_{2} \left(\frac{\left(\mu +\kappa \right)^{2} }{\sigma } +3\right)\left(\sum _{n=1}^{\infty }\frac{8}{n^{2} \pi ^{2} }  \right) \nonumber \\
&\times \left(\xi _{0}^{2} +w_{0}^{2} +\left\| \bar{\varphi }_{0} \right\| _{2}^{2} +\left\| \varphi _{0} \right\| _{H^{1} (0,1)}^{2} \right)
\end{align} 
Equations \eqref{GrindEQ__68_}, \eqref{GrindEQ__69_} and estimates \eqref{GrindEQ__72_}, \eqref{GrindEQ__73_}, \eqref{GrindEQ__74_} imply the existence of a constant $c_{3} >0$ (independent of the particular solution and $m$, $T>0$) for which the following estimate holds for all $t\in \left[0,T\right]$:
\begin{align} \label{GrindEQ__77_} 
\xi _{m}^{2} (t)+w_{m}^{2} (t)+f_{m}^{2} (t)+\left\| u_{m} [t]\right\| _{2}^{2} +\left\| \varphi _{m} [t]\right\| _{2}^{2} +\left\| \varphi _{m,x} [t]\right\| _{2}^{2} \nonumber \\ 
\le c_{3} \exp \left(-c_{1}^{-1} At\right)\left(\xi _{0}^{2} +w_{0}^{2} +\left\| \bar{\varphi }_{0} \right\| _{2}^{2} +\left\| \varphi _{0} \right\| _{H^{1} (0,1)}^{2} \right)
\end{align} 
Consequently, estimates \eqref{GrindEQ__76_}, \eqref{GrindEQ__77_} and the facts that $\left\| \varphi _{m} [t]\right\| _{2} \le \left\| \varphi _{m,x} [t]\right\| _{2} \le \left\| \varphi _{m,xx} [t]\right\| _{2} $ and $\left\| \varphi _{m,t} [t]\right\| _{2} \le \left\| \varphi _{m,tx} [t]\right\| _{2} $ for all $t\in \left[0,T\right]$, imply that there exists a constant $K>0$ (independent of $m$, $T>0$) such that the following inequality holds for all $t\in \left[0,T\right]$:
\begin{align} \label{GrindEQ__78_} 
\xi _{m}^{2} (t)+w_{m}^{2} (t)+f_{m}^{2} (t)+\left\| u_{m} [t]\right\| _{2}^{2} +\left\| \varphi _{m,t} [t]\right\| _{2}^{2} +\int _{0}^{t}f_{m}^{2} (s)ds \nonumber \\
 +\left\| \varphi _{m} [t]\right\| _{H^{2} (0,1)}^{2} +\int _{0}^{t}\left\| \varphi _{m} [s]\right\| _{H^{3} (0,1)}^{2} ds +\int _{0}^{t}\left\| \varphi _{m,t} [s]\right\| _{H^{1} (0,1)}^{2} ds \le K 
\end{align}
Equations \eqref{GrindEQ__28_}, \eqref{GrindEQ__14_}, \eqref{GrindEQ__15_}, \eqref{GrindEQ__62_}, \eqref{GrindEQ__29_} give for a.e. $t\in \left(0,T\right)$: 
\begin{equation} \label{GrindEQ__79_} 
u_{m,t} [t]=\varphi _{m,x} [t]-\kappa u_{m} [t]+\mu \varphi _{m,tx} [t]-\sigma \varphi _{m,xxx} [t]+f_{m} (t)\sum _{n=1}^{m}\frac{\beta _{n} }{n\pi } g_{n}   
\end{equation}
Equations \eqref{GrindEQ__14_}, \eqref{GrindEQ__28_}, \eqref{GrindEQ__30_} imply that $u_{m} (t,0)=u_{m} (t,1)=0$, $\left\| u_{m} [t]\right\| _{2} \le \left\| \varphi _{m,t} [t]\right\| _{2} $, $\left\| u_{m,x} [t]\right\| _{2} \le \left\| \varphi _{m,t} [t]\right\| _{2} $ and $\left\| u_{m,xx} [t]\right\| _{2} \le \left\| \varphi _{m,tx} [t]\right\| _{2} $ for all $t\in \left[0,T\right]$. These facts in conjunction with \eqref{GrindEQ__78_} and \eqref{GrindEQ__79_} implies that there exists a constant $\bar{K}>0$ (independent of $m$, $T>0$) such that the following inequality holds for all $t\in \left[0,T\right]$:
\begin{equation} \label{GrindEQ__80_} 
\left\| u_{m} [t]\right\| _{H_{0}^{1} (0,1)}^{2} +\int _{0}^{t}\left\| u_{m} [s]\right\| _{H^{2} (0,1)}^{2} ds +\int _{0}^{t}\left\| u_{m,t} [s]\right\| _{2}^{2} ds \le \bar{K} 
\end{equation} 
Working exactly as in the proof of Theorem \ref{theorem1} we establish inequalities \eqref{GrindEQ__33_}, \eqref{GrindEQ__34_}.

Thus, by virtue of \eqref{GrindEQ__78_} the sequence $\left\{\, \varphi _{m} \, :\, m=1,2,...\right\}$ is bounded in $L^{\infty } \left((0,T);H^{2} (0,1)\right)$ and in $L^{2} \left((0,T);H^{3} (0,1)\right)$. Moreover, the sequence 
$$\left\{\, \varphi _{m,t} \, :\, m=1,2,...\right\}$$ is bounded in $L^{\infty } \left((0,T);L^{2} (0,1)\right)$ and in $L^{2} \left((0,T);H^{1} (0,1)\right)$ while by virtue of \eqref{GrindEQ__80_} the sequence $\left\{\, u_{m} \, :\, m=1,2,...\right\}$ is bounded in $L^{\infty } \left((0,T);H_{0}^{1} (0,1)\right)$ and in $L^{2} \left((0,T);H^{2} (0,1)\right)$ and the sequence $\left\{\, u_{m,t} \, :\, m=1,2,...\right\}$ is bounded in $L^{2} \left((0,T);L^{2} (0,1)\right)$. By virtue of \eqref{GrindEQ__78_}, \eqref{GrindEQ__80_} and \eqref{GrindEQ__33_}, \eqref{GrindEQ__34_} the sequences $\left\{\, \varphi _{m} \, :\, m=1,2,...\right\}$, $\left\{\, u_{m} \, :\, m=1,2,...\right\}$ are uniformly bounded in $H^{1} (0,1)$ and are uniformly equicontinuous in $L^{2} (0,1)$. 

By virtue of \eqref{GrindEQ__78_} the sequences $\left\{\, f_{m} \, :\, m=1,2,...\right\}$, $\left\{\, \xi _{m} \, :\, m=1,2,...\right\}$ and $\left\{\, w_{m} \, :\, m=1,2,...\right\}$ are bounded in $C^{0} \left([0,T]\right)$. Furthermore, by virtue of \eqref{GrindEQ__63_} and \eqref{GrindEQ__66_} the sequences $\left\{\, \dot{f}_{m} \, :\, m=1,2,...\right\}$, $\left\{\, \dot{\xi }_{m} \, :\, m=1,2,...\right\}$ and $\left\{\, \dot{w}_{m} \, :\, m=1,2,...\right\}$ are bounded in $L^{2} (0,T)$ which also implies that the sequences $\left\{\, f_{m} \, :\, m=1,2,...\right\}$, $\left\{\, \xi _{m} \, :\, m=1,2,...\right\}$ and $\left\{\, w_{m} \, :\, m=1,2,...\right\}$ are uniformly equicontinuous. 

Consequently, there exist subsequences 
$$\left\{\, \varphi _{m_{k} } \, :\, k=1,2,...\right\}\subseteq \left\{\, \varphi _{m} \, :\, m=1,2,...\right\},$$ $$ \left\{\, u_{m_{k} } \, :\, k=1,2,...\right\}\subseteq \left\{\, u_{m} \, :\, m=1,2,...\right\},$$ $$ \left\{\, f_{m_{k} } \, :\, k=1,2,...\right\}\subseteq \left\{\, f_{m} \, :\, m=1,2,...\right\}, $$ $$ \left\{\, \xi _{m_{k} } \, :\, k=1,2,...\right\}\subseteq \left\{\, \xi _{m} \, :\, m=1,2,...\right\}, $$ $$ \left\{\, w_{m_{k} } \, :\, k=1,2,...\right\}\subseteq \left\{\, w_{m} \, :\, m=1,2,...\right\}$$
and 
$$\varphi \in L^{\infty } \left((0,T);H^{2} (0,1)\right)\cap L^{2} \left((0,T);H^{3} (0,1)\right), $$ $$ u\in L^{\infty } \left((0,T);H_{0}^{1} (0,1)\right)\cap L^{2} \left((0,T);H^{2} (0,1)\right), f,\xi ,w\in C^{0} \left([0,T]\right)$$ 
with 
$$\varphi _{t} \in L^{\infty } \left((0,T);L^{2} (0,1)\right)\cap L^{2} \left((0,T);H^{1} (0,1)\right),$$ $$u_{t} \in L^{2} \left((0,T);L^{2} (0,1)\right)$$ 
such that
\begin{align*}
\varphi _{m_{k} } \to \varphi  \textrm{ in } L^{2} \left((0,T);H^{3} (0,1)\right) \textrm{ weakly} \\
 u_{m_{k} } \to u \textrm{ in } L^{2} \left((0,T);H^{2} (0,1)\right) \textrm{ weakly}  \\
 \varphi _{m_{k} ,t} \to \varphi _{t}  \textrm{ in } L^{2} \left((0,T);H^{1} (0,1)\right) \textrm{ weakly} \\
u_{m_{k} ,t} \to u_{t}  \textrm{ in } L^{2} \left((0,T);L^{2} (0,1)\right) \textrm{ weakly} \\
\varphi _{m_{k} } \to \varphi \textrm{ and } u_{m_{k} } \to u \textrm{ in } C^{0} \left([0,T];L^{2} (0,1)\right) \textrm{ strongly} \\
\xi _{m_{k} } \to \xi , w_{m_{k} } \to w \textrm{ and } f_{m_{k} } \to f \textrm{ in } C^{0} \left([0,T]\right) \textrm{ strongly}
\end{align*}
The fact that $u\in C^{0} \left([0,T];H^{1} (0,1)\right)$ and $\varphi \in C^{0} \left([0,T];H^{2} (0,1)\right)$ follows from Theorem 4 on pages 287-288 in \cite{8}. Working exactly as in the proof of Theorem \ref{theorem1} we establish that $\varphi _{t} \in C^{0} \left([0,T];S\right)$, $\varphi \in C^{0} \left([0,T];\bar{S}\right)$, $u\in C^{0} \left([0,T];H_{0}^{1} (0,1)\right)$, $\varphi [0]=\varphi _{0} $, $\varphi _{t} [0]=\bar{\varphi }_{0} $ and that \eqref{GrindEQ__10_} holds. 

We show next that $\left(\xi (t),w(t)\right)$ for $t\in \left[0,T\right]$ and $\varphi $ constitute a weak solution on $[0,T]$ of the initial-boundary value problem \eqref{GrindEQ__1_}-\eqref{GrindEQ__5_} with initial condition \eqref{GrindEQ__9_} that satisfies \eqref{GrindEQ__58_} for $t\in [0,T]$.

First, we notice that since $\xi _{m} (t)=\xi _{0} +\int _{0}^{t}w_{m} (s)ds $, $w_{m} (t)=w_{0} -\int _{0}^{t}f_{m} (s)ds $ for all $m=1,2,...$, $t\in [0,T]$ and since $\xi _{m_{k} } \to \xi $, $w_{m_{k} } \to w$ and $f_{m_{k} } \to f$ in $C^{0} \left([0,T]\right)$ strongly we get that $\xi (t)=\xi _{0} +\int _{0}^{t}w(s)ds $, $w(t)=w_{0} -\int _{0}^{t}f(s)ds $ for all $t\in [0,T]$. Therefore, $\xi ,w\in C^{1} \left([0,T]\right)$ and \eqref{GrindEQ__1_} holds for all $t\in [0,T]$.

Exploiting \eqref{GrindEQ__30_} and \eqref{GrindEQ__79_} we get for every $v\in L^{2} \left((0,T);L^{2} (0,1)\right)$:
\begin{align} \label{GrindEQ__81_} 
&\int _{0}^{T}\int _{0}^{1}\left(u_{m,t} (t,x)-\varphi _{m,x} (t,x)+\kappa u_{m} (t,x)-\mu u_{m,xx} (t,x)\right)v(t,x)dxdt    \nonumber \\ 
&+\int _{0}^{T}\int _{0}^{1}\left(\sigma \varphi _{m,xxx} (t,x)+f(t)\right)v(t,x)dxdt \nonumber \\
&+\int _{0}^{T}\int _{0}^{1}\left(f_{m} (t)-f(t)\right)v(t,x)dxdt   \nonumber \\
&-\int _{0}^{T}\int _{0}^{1}f_{m} (t)\left(1+\sqrt{2} \sum _{n=1}^{m}\frac{(-1)^{n} -1}{n\pi } g_{n} (x) \right)v(t,x)dxdt  =0
\end{align}
Taking into account the fact that the functions $g_{n} $ defined by \eqref{GrindEQ__29_} constitute an orthonormal basis of $L^{2} (0,1)$ and Parseval's identity, we get 
\[\left\| \chi _{[0,1]} -\sum _{n=1}^{m}\frac{\left(\phi _{n} (0)-\phi _{n} (1)\right)}{n\pi } g_{n}  \right\| _{2}^{2} =\sum _{n=m+1}^{\infty }\frac{\left(\phi _{n} (0)-\phi _{n} (1)\right)^{2} }{n^{2} \pi ^{2} }  \le 4\sum _{n=m+1}^{\infty }\frac{1}{n^{2} \pi ^{2} }  \] 
which combined with \eqref{GrindEQ__15_} and the Cauchy-Schwarz inequality gives: 
\begin{align} \label{GrindEQ__82_} 
&\left|\int _{0}^{T}\int _{0}^{1}\left(f_{m}(t)-f(t) \right)v(t,x)dxdt  \right. \nonumber \\
&\left. -\int _{0}^{T}\int _{0}^{1}f_{m} (t)\left(1+\sqrt{2} \sum _{n=1}^{m}\frac{(-1)^{n} -1}{n\pi } g_{n} (x) \right)v(t,x)dxdt  \right| \nonumber \\ 
\le& 2\left(\sum _{n=m+1}^{\infty }\frac{1}{n^{2} \pi ^{2} }  \right)^{1/2} \left(\int _{0}^{T}\left|f_{m} (t)\right|^{2} dt \right)^{1/2} \left(\int _{0}^{T}\left\| v[t]\right\| _{2}^{2} dt \right)^{1/2} \nonumber \\ 
&+\left(\int _{0}^{T}\left|f_{m} (t)-f(t)\right|^{2} dt \right)^{1/2} \left(\int _{0}^{T}\left\| v[t]\right\| _{2}^{2} dt \right)^{1/2}   
\end{align}
Therefore, it follows from \eqref{GrindEQ__78_}, \eqref{GrindEQ__82_} and the fact that $f_{m_{k} } \to f$ in $C^{0} \left([0,T]\right)$ strongly that
\begin{align}
&{\mathop{\lim }\limits_{k\to +\infty }} \left(\int _{0}^{T}\int _{0}^{1}\left(f_{m_{k}} (t)-f(t) \right)v(t,x)dxdt  \right. \nonumber \\
&\left. -\int _{0}^{T}\int _{0}^{1}f_{m_{k} } (t)\left(1+\sqrt{2} \sum _{n=1}^{m_{k} }\frac{(-1)^{n} -1}{n\pi } g_{n} (x) \right)v(t,x)dxdt  \right)=0
\end{align}
Since $\varphi _{m_{k} } \to \varphi $ in $L^{2} \left((0,T);H^{3} (0,1)\right)$ weakly and $u_{m_{k} } \to u$ in $L^{2} \left((0,T);\right.$ $\left. H^{2} (0,1)\right)$ weakly, we get from \eqref{GrindEQ__81_}:
\begin{align} \label{GrindEQ__83_} 
&\int _{0}^{T}\int _{0}^{1}\left(u_{t} (t,x)-\varphi _{x} (t,x)+\kappa u(t,x)-\mu u_{xx} (t,x)\right)v(t,x)dxdt =0 \nonumber \\
&+\int _{0}^{T}\int _{0}^{1}\left(\sigma \varphi _{xxx} (t,x)+f(t)\right)v(t,x)dxdt =0
\end{align} 
Notice that $u_{t} [t]-\varphi _{x} [t]+\kappa u[t]-\mu u_{xx} [t]+\sigma \varphi _{xxx} [t]+f(t)\chi _{[0,1]} $ is of class $L^{2} \left((0,T);L^{2} (0,1)\right)$. Since \eqref{GrindEQ__83_} holds for every $v\in L^{2} \left((0,T);L^{2} (0,1)\right)$, we obtain \eqref{GrindEQ__11_}.

Moreover, \eqref{GrindEQ__66_} implies that for all $m=1,2,...$ the following equation holds for all $v\in L^{2} \left(0,T\right)$:
\begin{align} \label{GrindEQ__84_} 
&\int _{0}^{T}v(t)\left(\gamma ^{-1} k^{-3} f_{m} (t)-k^{-2} \left(w_{m} (t)+k\xi _{m} (t)\right)\right)dt  \nonumber \\ 
&=-rQ\int _{0}^{T}\int _{0}^{1}v(t)\left(\kappa x\varphi _{m} (t,x)-\frac{r+1}{r} u_{m} (t,x)+\mu \varphi _{m,x} (t,x)\right)dx dt  
\end{align} 
Since $\varphi _{m_{k} } \to \varphi $ in $L^{2} \left((0,T);H^{3} (0,1)\right)$ weakly, $u_{m_{k} } \to u$ in $L^{2} \left((0,T);H^{2} (0,1)\right)$ weakly and $\xi _{m_{k} } \to \xi $, $w_{m_{k} } \to w$, $f_{m_{k} } \to f$ in $C^{0} \left([0,T];{\mathbb R}\right)$ strongly, we obtain from \eqref{GrindEQ__84_} that the following equation holds for all $v\in L^{2} \left(0,T\right)$:
\begin{align} \label{GrindEQ__85_} 
&\int _{0}^{T}v(t)\left(\gamma ^{-1} k^{-3} f(t)-k^{-2} \left(w(t)+k\xi (t)\right)\right)dt  \nonumber \\ 
&=-rQ\int _{0}^{T}\int _{0}^{1}v(t)\left(\kappa x\varphi (t,x)-\frac{r+1}{r} u(t,x)+\mu \varphi _{x} (t,x)\right)dx dt 
\end{align}
 Consequently, we obtain the following equation for all $t\in \left[0,T\right]$ (continuity of $\varphi \in C^{0} \left([0,T];\bar{S}\right)$ and $u\in C^{0} \left([0,T];H_{0}^{1} (0,1)\right)$ is also exploited at this point):
\begin{align} \label{GrindEQ__86_} 
f(t)=& \gamma k \left(w(t)+k\xi (t)\right) \nonumber \\
& -\gamma k^{3} rQ\int _{0}^{1}\left(\kappa x\varphi (t,x)-\frac{r+1}{r} u(t,x)+\mu \varphi _{x} (t,x)\right)dx 
\end{align}
Equation \eqref{GrindEQ__58_} for $t\in [0,T]$ is a direct consequence of \eqref{GrindEQ__86_}, \eqref{GrindEQ__56_} and \eqref{GrindEQ__10_}.

We next show that every weak solution on $[0,T]$ of the initial-boundary value problem \eqref{GrindEQ__1_}-\eqref{GrindEQ__5_} with initial condition \eqref{GrindEQ__9_} that satisfies \eqref{GrindEQ__58_} for $t\in [0,T]$ also satisfies \eqref{GrindEQ__59_} and \eqref{GrindEQ__60_}. Uniqueness of weak solution follows directly from \eqref{GrindEQ__59_}.

The mapping $\left[0,T\right]\backepsilon t\mapsto V\left(\xi (t),w(t),\varphi [t],u[t]\right)$ with $V$ being defined by \eqref{GrindEQ__57_} is absolutely continuous and satisfies the following equation for a.e. $t\in \left(0,T\right)$:  
\begin{align} \label{GrindEQ__87_} 
&\frac{d}{d\, t} V\left(\xi (t),w(t),\varphi [t],u[t]\right) \nonumber \\
=&\xi (t)\dot{\xi }(t)+k^{-2} \left(w(t)+k\xi (t)\right)\left(\dot{w}(t)+k\dot{\xi }(t)\right) \nonumber \\ 
&+Q\left(u[t],u_{t} [t]\right)+Q(1+r)\left(\varphi [t],\varphi _{t} [t]\right)+Q\sigma (1+r)\left(\varphi _{x} [t],\varphi _{tx} [t]\right) \nonumber \\ 
&+Qr\left(u[t]-\mu \varphi _{x} [t]+\kappa \mho \varphi [t],u_{t} [t]-\mu \varphi _{tx} [t]+\kappa \left(\mho \varphi [t]\right)_{t} \right) 
\end{align} 
Using \eqref{GrindEQ__1_}, \eqref{GrindEQ__8_}, \eqref{GrindEQ__10_}, \eqref{GrindEQ__11_} (which also imply the equalities $\left(\mho \varphi [t]\right)_{t} =u[t]$, $\varphi _{tx} [t]=u_{xx} [t]$, $\left(\mho \varphi [t]\right)_{x} =\varphi [t]$) and integration by part we obtain from \eqref{GrindEQ__87_} for a.e. $t\in \left(0,T\right)$:
\begin{align}  \label{GrindEQ__88_}
&\frac{d}{d\, t} V\left(\xi (t),w(t),\varphi [t],u[t]\right) \nonumber \\
=&\xi (t)w(t)+k^{-1} \left(w(t)+k\xi (t)\right)w(t)-Q\kappa \left\| u[t]\right\| _{2}^{2} -Q\mu \left\| u_{x} [t]\right\| _{2}^{2} \nonumber  \\ 
&-Q\kappa r\left\| \varphi [t]\right\| _{2}^{2} -Qr\left(\sigma \kappa +\mu \right)\left\| \varphi _{x} [t]\right\| _{2}^{2} -Q\mu \sigma r\left\| \varphi _{xx} [t]\right\| _{2}^{2} \nonumber \\ 
&-f(t)\left(k^{-2} \left(w(t)+k\xi (t)\right)-Qr\kappa \int _{0}^{1} x\varphi (t,x) dx \right) \nonumber \\
&-f(t)\left(Q(r+1)\int _{0}^{1}u(t,x)dx -Qr\mu \left(\varphi (t,1)-\varphi (t,0)\right)\right)
\end{align}
Using 
\eqref{GrindEQ__86_} and \eqref{GrindEQ__88_} we get for a.e. $t\in \left(0,T\right)$:
\begin{align} \label{GrindEQ__89_} 
&\frac{d}{d\, t} V\left(\xi (t),w(t),\varphi [t],u[t]\right) = \nonumber \\
&-k\xi ^{2} (t)+k^{-1} \left(w(t)+k\xi (t)\right)^{2} -Q\kappa \left\| u[t]\right\| _{2}^{2} -Q\mu \left\| u_{x} [t]\right\| _{2}^{2}   \nonumber \\ 
&-Q\kappa r\left\| \varphi [t]\right\| _{2}^{2} -Qr\left(\sigma \kappa +\mu \right)\left\| \varphi _{x} [t]\right\| _{2}^{2} -Q\mu \sigma r\left\| \varphi _{xx} [t]\right\| _{2}^{2} -\frac{1}{2\gamma k^{3} } f^{2} (t)  \nonumber \\ 
&-\hspace{-0.2em} \frac{\gamma k^{3} }{2} \hspace{-0.3em} \left( \hspace{-0.2em} k^{-2}  \left(w(t)+k\xi (t)\right)-Qr\int _{0}^{1}\hspace{-0.4em} \left(\kappa x\varphi (t,x)-\frac{r+1}{r} u(t,x)+\mu \varphi _{x} (t,x)\right) \hspace{-0.2em} dx \hspace{-0.1em} \right) ^{2}   
\end{align} 
Definition \eqref{GrindEQ__40_} in conjunction with \eqref{GrindEQ__89_} implies the equality
\begin{equation*}
\frac{d}{d\, t} V\left(\xi (t),w(t),\varphi [t],u[t]\right)+q\left(\xi (t),w(t),\varphi [t],u[t]\right)+\frac{1}{2\gamma k^{3} } f^{2} (t)=0 
\end{equation*} 
for a.e. $t\in \left(0,T\right)$. Integrating this equation and using \eqref{GrindEQ__8_}, \eqref{GrindEQ__9_}, \eqref{GrindEQ__10_}, we obtain \eqref{GrindEQ__60_}. Definitions \eqref{GrindEQ__12_}, \eqref{GrindEQ__57_} as well as the fact $\pi ^{2} \left\| \mho \varphi \right\| _{2}^{2} \le \left\| \varphi \right\| _{2}^{2} $ (Wirtinger's inequality; recall definition \eqref{GrindEQ__8_}) imply the existence of a constant $c_{5} >0$ for which the following inequality holds for all $\varphi \in \bar{S}$, $u\in H_{0}^{1} (0,1)$, $\left(\xi ,w\right)\in {\mathbb R}^{2} $: 
\begin{equation} \label{GrindEQ__90_} 
V\left(\xi ,w,\varphi ,u\right)\le c_{5} \left(\xi ^{2} +w^{2} +\left\| u\right\| _{2}^{2} +\left\| \varphi \right\| _{H^{1} (0,1)}^{2} \right) 
\end{equation}
Exploiting \eqref{GrindEQ__B1_}, \eqref{GrindEQ__58_},   \eqref{GrindEQ__60_} and \eqref{GrindEQ__90_} we conclude that the differential inequality 
\begin{equation*}
\frac{d}{d\, t} V\left(\xi (t),w(t),\varphi [t],u[t]\right)\le -c_{5}^{-1} \bar{A}V\left(\xi (t),w(t),\varphi [t],u[t]\right) 
\end{equation*}
holds for all  $t\in \left[0,T\right]$. Consequently, we get estimate \eqref{GrindEQ__59_} with $c=c_{5}^{-1} \bar{A}$. The proof is complete.
\qed
\end{proof}

Next, we give the proof of Theorem \ref{theorem3}. 

\begin{proof}{\hspace{-0.15em}\textit{of Theorem \ref{theorem3}}}
 Let arbitrary $\varphi _{0} \in \bar{S}$, $\bar{\varphi }_{0} \in S$, $\left(\xi _{0} ,w_{0} \right)\in {\mathbb R}^{2} $ be given. Using \eqref{GrindEQ__13_}, \eqref{GrindEQ__57_} and \eqref{GrindEQ__1_} we conclude that for every $f\in L_{loc}^{2} $ the corresponding unique weak solution on ${\mathbb R}_{+} $ of the initial-boundary value problem \eqref{GrindEQ__1_}-\eqref{GrindEQ__5_} with initial condition \eqref{GrindEQ__9_} satisfies the following equation for all $t\ge 0$:
 \begin{align} \label{GrindEQ__96_} 
V\left(\xi (t),w(t),\varphi [t],u[t]\right) 
=& V\left(\xi _{0} ,w_{0} ,\varphi _{0} ,\mho \bar{\varphi }_{0} \right) \nonumber \\
&+\int _{0}^{t}\xi (s)w(s)ds +k^{-1} \int _{0}^{t}w(s)\left(w(s)+k\xi (s)\right)ds \nonumber  \\ 
&-Q\int _{0}^{t}\left(r\left(\mu +\kappa \sigma \right)\left\| \varphi _{x} [s]\right\| _{2}^{2} +r\mu \sigma \left\| \varphi _{xx} [s]\right\| _{2}^{2}  \right)ds \nonumber \\
&-Q\int _{0}^{t}\left( \kappa r\left\| \varphi [s]\right\| _{2}^{2} +\kappa \left\| u[s]\right\| _{2}^{2} +\mu \left\| u_{x} [s]\right\| _{2}^{2} \right)ds  \nonumber \\ 
&-\int _{0}^{t}f(s)\left(k^{-2} \left(w(s)+k\xi (s)\right) \right) ds   \nonumber \\  
&+Qr\int _{0}^{t}f(s)\left(\int _{0}^{1} \kappa x\varphi (s,x) dx \right)ds  \nonumber \\
&-Qr\int _{0}^{t}f(s)\left(\int _{0}^{1} \frac{r+1}{r} u(s,x) dx \right)ds   \nonumber \\
&+Qr \int _{0}^{t}f(s)\left(\int _{0}^{1} \mu \varphi _{x} (s,x) dx \right)ds 
\end{align} 
where $u\in C^{0} \left([0,T];H_{0}^{1} (0,1)\right)\cap L^{2} \left((0,T);H^{2} (0,1)\right)$ is the function that satisfies \eqref{GrindEQ__10_}, \eqref{GrindEQ__11_}. Using \eqref{GrindEQ__40_}, \eqref{GrindEQ__56_} and \eqref{GrindEQ__10_} (which implies that $\int _{0}^{1}u(s,x)dx =-\int _{0}^{1}x\varphi _{t} (s,x)dx $) we get from \eqref{GrindEQ__96_} the following equation for all $t\ge 0$:  
\begin{align} \label{GrindEQ__97_} 
&V\left(\xi (t),w(t),\varphi [t],u[t]\right)+\int _{0}^{t}q\left(\xi (s),w(s),\varphi [s],u[s]\right)ds +\frac{1}{2\gamma k^{3} } \int _{0}^{t}f^{2} (s)ds \nonumber \\ 
&=V\left(\xi _{0} ,w_{0} ,\varphi _{0} ,\mho \bar{\varphi }_{0} \right)+\frac{1}{2\gamma k^{3} } \int _{0}^{t}\left(f(s)-P\left(\xi (s),w(s),\varphi [s],\varphi _{t} [s]\right)\right)^{2} ds  
\end{align} 
We next claim that for every $f\in L_{loc}^{2} $ the following inequality holds: 
\begin{equation} \label{GrindEQ__98_} 
J\left(\xi _{0} ,w_{0} ,\varphi _{0} ,\bar{\varphi }_{0} ,f\right)\ge V\left(\xi _{0} ,w_{0} ,\varphi _{0} ,\mho \bar{\varphi }_{0} \right) 
\end{equation} 
Indeed, \eqref{GrindEQ__98_} holds automatically when $J\left(\xi _{0} ,w_{0} ,\varphi _{0} ,\bar{\varphi }_{0} ,f\right)=+\infty $. So, we need to prove \eqref{GrindEQ__98_} only for the case where $J\left(\xi _{0} ,w_{0} ,\varphi _{0} ,\bar{\varphi }_{0} ,f\right)<+\infty $. So, next we assume that $J\left(\xi _{0} ,w_{0} ,\varphi _{0} ,\bar{\varphi }_{0} ,f\right)<+\infty $. Equation \eqref{GrindEQ__97_} gives for $t\ge 0$ a.e.:
\begin{align} \label{GrindEQ__99_} 
\frac{d}{d\, t} V\left(\xi (t),w(t),\varphi [t],u[t]\right)=&-q\left(\xi (t),w(t),\varphi [t],u[t]\right) \nonumber \\ 
&+\frac{1}{2\gamma k^{3} } P^{2} \left(\xi (t),w(t),\varphi [t],\varphi _{t} [t]\right) \nonumber \\
&-\frac{1}{\gamma k^{3} } f(t)P\left(\xi (t),w(t),\varphi [t],\varphi _{t} [t]\right)  
\end{align} 
Definitions \eqref{GrindEQ__12_}, \eqref{GrindEQ__57_} as well as the fact $\pi ^{2} \left\| \mho \varphi \right\| _{2}^{2} \le \left\| \varphi \right\| _{2}^{2} $ (Wirtinger's inequality; recall definition \eqref{GrindEQ__8_}) imply the existence of a constant $c_{5} >0$ for which estimate \eqref{GrindEQ__90_} holds for all $\varphi \in \bar{S}$, $u\in H_{0}^{1} (0,1)$, $\left(\xi ,w\right)\in {\mathbb R}^{2} $. Definition \eqref{GrindEQ__56_} implies the existence of a constant $c_{6} >0$ for which the following estimate holds for all $\varphi \in \bar{S}$, $\bar{\varphi }\in S$, $\left(\xi ,w\right)\in {\mathbb R}^{2} $:
\begin{equation} \label{GrindEQ__100_} 
\left|P\left(\xi ,w,\varphi ,\bar{\varphi }\right)\right|^{2} \le c_{6} \left(\xi ^{2} +w^{2} +\left\| \bar{\varphi }\right\| _{2}^{2} +\left\| \varphi \right\| _{H^{1} (0,1)}^{2} \right) 
\end{equation} 
Using the inequality 
$$-f(t)P\left(\xi (t),w(t),\varphi [t],\varphi _{t} [t]\right)\le \frac{1}{2} P^{2} \left(\xi (t),w(t),\varphi [t],\varphi _{t} [t]\right)+\frac{1}{2} f^{2} (t)$$
in conjunction with inequalities \eqref{GrindEQ__42_}, \eqref{GrindEQ__90_}, \eqref{GrindEQ__100_} and the fact that $\left\| u[t]\right\| _{H^{1} (0,1)}^{2} \ge \left\| \varphi _{t} [t]\right\| _{2}^{2} $ (recall equation \eqref{GrindEQ__10_}), we obtain from \eqref{GrindEQ__99_} for $t\ge 0$ a.e.:
\begin{align} \label{GrindEQ__101_} 
\frac{d}{d\, t} V\left(\xi (t),w(t),\varphi [t],u[t]\right)\le & -c_{5}^{-1} AV\left(\xi (t),w(t),\varphi [t],u[t]\right) \nonumber \\ 
&+\left(\frac{c_{6} }{A\gamma k^{3}} +1\right)q\left(\xi (t),w(t),\varphi [t],u[t]\right) \nonumber \\
&+\frac{1}{2\gamma k^{3} } \left(\frac{c_{6} }{A\gamma k^{3}} +1\right) f^{2} (t)
\end{align} 
Integrating the differential inequality \eqref{GrindEQ__101_} and using definition \eqref{GrindEQ__91_} as well as \eqref{GrindEQ__9_} and \eqref{GrindEQ__10_}, we get for all $t\ge T\ge 0$:  
\begin{align} \label{GrindEQ__102_} 
& V\left(\xi (t),w(t),\varphi [t],u[t]\right) \nonumber \\
\le & \exp \left(-c_{5}^{-1} At\right)V\left(\xi _{0} ,w_{0} ,\varphi _{0} ,\mho \bar{\varphi }_{0} \right) \nonumber \\ 
&+\left(\frac{c_{6} }{A\gamma k^{3}} +1\right)\int _{0}^{t}\exp \left(-c_{5}^{-1} A(t-s)\right)q\left(\xi (s),w(s),\varphi [s],u[s]\right)ds  \nonumber \\
&+\frac{1}{2\gamma k^{3} } \left(\frac{c_{6} }{A\gamma k^{3}} +1\right)  \int _{0}^{t}\exp \left(-c_{5}^{-1} A(t-s)\right)f^{2} (s) ds  \nonumber \\
\le& \exp \left(-c_{5}^{-1} At\right)V\left(\xi _{0} ,w_{0} ,\varphi _{0} ,\mho \bar{\varphi }_{0} \right) \nonumber \\ 
&+\exp \left(-c_{5}^{-1} A(t-T)\right)\left(\frac{c_{6} }{A\gamma k^{3} } +1\right)J\left(\xi _{0} ,w_{0} ,\varphi _{0} ,\bar{\varphi }_{0} ,f\right) \nonumber \\ 
&+\left(\frac{c_{6} }{A\gamma k^{3} } +1\right)\int _{T}^{t}\left(q\left(\xi (s),w(s),\varphi [s],u[s]\right)+\frac{1}{2\gamma k^{3}} f^{2} (s)\right)ds   
\end{align} 
Estimate \eqref{GrindEQ__102_} shows that
\begin{align*}
&{\mathop{\lim \sup }\limits_{t\to +\infty }} \left(V\left(\xi (t),w(t),\varphi [t],u[t]\right)\right) \nonumber \\
&\le \left(\frac{c_{6} }{A\gamma k^{3}} +1\right)\int _{T}^{+\infty }\left(q\left(\xi (s),w(s),\varphi [s],u[s]\right)+\frac{1}{2\gamma k^{3}} f^{2} (s)\right)ds 
\end{align*}
for all $T\ge 0$. Since $\int _{0}^{+\infty }\left(q\left(\xi (s),w(s),\varphi [s],u[s]\right)+\frac{1}{2\gamma k^{3} } f^{2} (s)\right)ds <+\infty $ (recall \eqref{GrindEQ__91_}) the above estimate shows that ${\mathop{\lim }\limits_{t\to +\infty }} \left(V\left(\xi (t),w(t),\varphi [t],u[t]\right)\right)=0$. Since \eqref{GrindEQ__97_} in conjunction with \eqref{GrindEQ__91_} gives the inequality $J\left(\xi _{0} ,w_{0} ,\varphi _{0} ,\mho \bar{\varphi }_{0} ,f\right)\ge V\left(\xi _{0} ,w_{0} ,\varphi _{0} ,\mho \bar{\varphi }_{0} \right)-V\left(\xi (t),w(t),\varphi [t],u[t]\right)$ for all $t\ge 0$, we obtain \eqref{GrindEQ__98_} by letting $t\to +\infty $. 
Equation \eqref{GrindEQ__97_} in conjunction with \eqref{GrindEQ__59_}, definition \eqref{GrindEQ__91_} and \eqref{GrindEQ__98_} implies that equation \eqref{GrindEQ__92_} holds. The proof is complete.
\qed
\end{proof}
Finally, we provide the proof of Theorem \ref{theorem4}.
\begin{proof}{\hspace{-0.15em}\textit{of Theorem \ref{theorem4}}} 
Clearly, by virtue of \eqref{GrindEQ__5_}, \eqref{GrindEQ__10_}, \eqref{GrindEQ__15_}, \eqref{GrindEQ__29_} and Parseval's identity, the following equations hold for $t\ge 0$: 
\begin{align}  \label{GrindEQ__104_}
\left\| \varphi _{t} [t]\right\| _{2}^{2} &=\left\| u_{x} [t]\right\| _{2}^{2} =\sum _{n=1}^{\infty }n^{2} \pi ^{2} \left(u[t],g_{n} \right)^{2}  \nonumber \\ 
\left\| \varphi [t]\right\| _{2}^{2} &=\sum _{n=1}^{\infty }\left(\varphi [t],\phi _{n} \right)^{2}    \nonumber \\  
\left\| \varphi _{x} [t]\right\| _{2}^{2} &=\sum _{n=1}^{\infty }n^{2} \pi ^{2} \left(\varphi [t],\phi _{n} \right)^{2}    \nonumber  \\
 \left\| \varphi _{xx} [t]\right\| _{2}^{2} &=\sum _{n=1}^{\infty }n^{4} \pi ^{4} \left(\varphi [t],\phi _{n} \right)^{2}  
\end{align} 
where $u\in C^{0} \left({\mathbb R}_{+} ;H_{0}^{1} (0,1)\right)$ is the function that satisfies \eqref{GrindEQ__10_}, \eqref{GrindEQ__11_} for $t\ge 0$. Exploiting \eqref{GrindEQ__10_}, \eqref{GrindEQ__11_}, \eqref{GrindEQ__15_}, \eqref{GrindEQ__29_} and \eqref{GrindEQ__65_} and using integration by parts, we obtain the following equations hold for $t\ge 0$ a.e.:
\begin{align} \label{GrindEQ__105_} 
\frac{d}{d\, t} \left(\varphi [t],\phi _{n} \right)&=n\pi \left(u[t],g_{n} \right) \nonumber \\ 
\frac{d}{d\, t} \left(u[t],g_{n} \right)&=-n\pi \left(1+\sigma n^{2} \pi ^{2} \right)\left(\varphi [t],\phi _{n} \right)-\left(\mu n^{2} \pi ^{2} +\kappa \right)\left(u[t],g_{n} \right)+\frac{\beta _{n} }{n\pi } f(t) 
\end{align} 
Define for $t\ge 0$ and $n=1,2,...$
\begin{align} \label{GrindEQ__106_} 
P_{n} (t):=&\frac{1}{2} n^{2} \pi ^{2} \left(u[t],g_{n} \right)^{2} +n^{2} \pi ^{2} \left(1+\sigma n^{2} \pi ^{2} \right)\left(\varphi [t],\phi _{n} \right)^{2}  \nonumber \\ 
&+\frac{1}{2} n^{2} \pi ^{2} \left(\left(u[t],g_{n} \right)+\frac{\mu n^{2} \pi ^{2} +\kappa }{n\pi } \left(\varphi [t],\phi _{n} \right)\right)^{2}  
\end{align}
Using \eqref{GrindEQ__105_} and \eqref{GrindEQ__106_} we get for $t\ge 0$ a.e.:
\begin{align} \label{GrindEQ__107_} 
 \dot{P}_{n} (t)=&-\left(\mu n^{2} \pi ^{2} +\kappa \right)n^{2} \pi ^{2} \left(u[t],g_{n} \right)^{2} \nonumber  \\
&-\left(1+\sigma n^{2} \pi ^{2} \right)\left(\mu n^{2} \pi ^{2} +\kappa \right)n^{2} \pi ^{2} \left(\varphi [t],\phi _{n} \right)^{2}    \nonumber \\
&+n\pi \beta _{n} \left(2\left(u[t],g_{n} \right)+\frac{\mu n^{2} \pi ^{2} +\kappa }{n\pi } \left(\varphi [t],\phi _{n} \right)\right)f(t)
\end{align}
Completing the squares in \eqref{GrindEQ__107_} we obtain the following inequality:
\begin{align} \label{GrindEQ__108_} 
 \dot{P}_{n} (t)\le & -\frac{n^{2} \pi ^{2} }{2} \left(\mu n^{2} \pi ^{2} +\kappa \right)\left(\left(u[t],g_{n} \right)^{2} +\left(1+\sigma n^{2} \pi ^{2} \right)\left(\varphi [t],\phi _{n} \right)^{2} \right) \nonumber \\ 
&+\left(\frac{2}{\mu n^{2} \pi ^{2} +\kappa } +\frac{\mu n^{2} \pi ^{2} +\kappa }{2\left(1+\sigma n^{2} \pi ^{2} \right)n^{2} \pi ^{2} } \right)\beta _{n}^{2} f^{2} (t) 
\end{align} 
Definition \eqref{GrindEQ__106_} implies the following inequality:
\begin{equation} \label{GrindEQ__109_} 
P_{n} (t)\le Cn^{2} \pi ^{2} \left(\left(u[t],g_{n} \right)^{2} +\left(1+\sigma n^{2} \pi ^{2} \right)\left(\varphi [t],\phi _{n} \right)^{2} \right) 
\end{equation} 
where 
\begin{equation} \label{GrindEQ__110_} 
C:=\frac{\left(\mu +\kappa \right)^{2} }{\sigma } +\frac{3}{2}  
\end{equation} 
It follows from \eqref{GrindEQ__108_} and \eqref{GrindEQ__109_} that the following differential inequality holds for $t\ge 0$ a.e.:
\begin{equation} \label{GrindEQ__111_} 
\dot{P}_{n} (t)\le -\frac{\mu n^{2} \pi ^{2} +\kappa }{2C} P_{n} (t)+\left(\frac{2}{\mu n^{2} \pi ^{2} +\kappa } +\frac{\mu n^{2} \pi ^{2} +\kappa }{2\left(1+\sigma n^{2} \pi ^{2} \right)n^{2} \pi ^{2} } \right)\beta _{n}^{2} f^{2} (t) 
\end{equation} 
Integrating \eqref{GrindEQ__111_} we obtain the following estimate for $t\ge 0$:
\begin{align} \label{GrindEQ__112_} 
P_{n} (t)\le& \exp \left(-\frac{\mu n^{2} \pi ^{2} +\kappa }{2C} t\right)P_{n} (0)  \nonumber \\
&+\left(\frac{2}{\mu } +\frac{\mu +\kappa }{2\sigma } \right)\frac{\beta _{n}^{2} }{n^{2} \pi ^{2} } \int _{0}^{t}\exp \left(-\frac{\mu n^{2} \pi ^{2} +\kappa }{2C} (t-s)\right)f^{2} (s)ds  
\end{align}
Using \eqref{GrindEQ__10_}, \eqref{GrindEQ__56_}, \eqref{GrindEQ__58_} and integration by parts we obtain \eqref{GrindEQ__86_} for $t\ge 0$. Exploiting \eqref{GrindEQ__12_}, \eqref{GrindEQ__57_}, \eqref{GrindEQ__58_} and \eqref{GrindEQ__86_} we conclude that there exists a constant $c_{7} >0$ (independent of the particular solution) such that the following estimate holds for $t\ge 0$: 
\begin{equation} \label{GrindEQ__113_} 
f^{2} (t)\le c_{7} V\left(\xi (t),w(t),\varphi [t],u[t]\right) 
\end{equation} 
Combining \eqref{GrindEQ__65_}, \eqref{GrindEQ__112_}, \eqref{GrindEQ__113_} and \eqref{GrindEQ__59_} we obtain the following estimate for $t\ge 0$:
\begin{align} \label{GrindEQ__114_} 
P_{n} (t)\le & \exp \left(-\frac{\mu \pi ^{2} +\kappa }{2C} t\right)P_{n} (0) \nonumber \\ 
&+\left(\frac{2}{\mu } +\frac{\mu +\kappa }{2\sigma } \right)\frac{8c_{7} }{n^{2} \pi ^{2} } \left(\int _{0}^{t}\exp \left(-\frac{\mu \pi ^{2} +\kappa }{2C} (t-s)\right)\exp \left(-cs\right)ds \right) \nonumber \\
&\times V\left(\xi _{0} ,w_{0} ,\varphi _{0} ,\mho \bar{\varphi }_{0} \right)
\end{align} 
Define $\displaystyle{\omega =\frac{1}{2} \min \left(c,\frac{\mu \pi ^{2} +\kappa }{2C} \right)}$. Using \eqref{GrindEQ__9_}, \eqref{GrindEQ__10_}, \eqref{GrindEQ__114_}, \eqref{GrindEQ__106_} and \eqref{GrindEQ__109_} we get for $t\ge 0$:
\begin{align} \label{GrindEQ__115_} 
&\frac{1}{2} n^{2} \pi ^{2} \left(u[t],g_{n} \right)^{2} +n^{2} \pi ^{2} \left(1+\sigma n^{2} \pi ^{2} \right)\left(\varphi [t],\phi _{n} \right)^{2} \nonumber  \\ 
\le& C\exp \left(-2\omega t\right)n^{2} \pi ^{2} \left(\left(\mho \bar{\varphi }_{0} ,g_{n} \right)^{2} +\left(1+\sigma n^{2} \pi ^{2} \right)\left(\varphi _{0} ,\phi _{n} \right)^{2} \right) \nonumber \\ 
&+\left(\frac{2}{\mu } +\frac{\mu +\kappa }{2\sigma } \right)\frac{8c_{7} }{n^{2} \pi ^{2} } t\exp \left(-2\omega t\right)V\left(\xi _{0} ,w_{0} ,\varphi _{0} ,\mho \bar{\varphi }_{0} \right)  
\end{align} 
Using the fact that $t\exp \left(-\omega t\right)\le \frac{1}{e\omega } $ for all $t\ge 0$, we get from \eqref{GrindEQ__115_} for $t\ge 0$:
\begin{align} \label{GrindEQ__116_} 
&\frac{1}{2} n^{2} \pi ^{2} \left(u[t],g_{n} \right)^{2} +n^{2} \pi ^{2} \left(1+\sigma n^{2} \pi ^{2} \right)\left(\varphi [t],\phi _{n} \right)^{2} \nonumber \\ 
\le& C\exp \left(-\omega t\right)n^{2} \pi ^{2} \left(\left(\mho \bar{\varphi }_{0} ,g_{n} \right)^{2} +\left(1+\sigma n^{2} \pi ^{2} \right)\left(\varphi _{0} ,\phi _{n} \right)^{2} \right) \nonumber \\ 
&+\left(\frac{2}{\mu } +\frac{\mu +\kappa }{2\sigma } \right)\frac{8c_{7} }{e\omega n^{2} \pi ^{2} } \exp \left(-\omega t\right)V\left(\xi _{0} ,w_{0} ,\varphi _{0} ,\mho \bar{\varphi }_{0} \right)
\end{align} 
Definition \eqref{GrindEQ__8_} and \eqref{GrindEQ__15_}, \eqref{GrindEQ__29_} imply that $n^{2} \pi ^{2} \left(\mho \bar{\varphi }_{0} ,g_{n} \right)^{2} =\left(\bar{\varphi }_{0} ,\phi _{n} \right)^{2} $. Definitions \eqref{GrindEQ__12_}, \eqref{GrindEQ__57_} as well as the fact $\pi ^{2} \left\| \mho \varphi \right\| _{2}^{2} \le \left\| \varphi \right\| _{2}^{2} $ (Wirtinger's inequality; recall definition \eqref{GrindEQ__8_}) imply the existence of a constant $c_{5} >0$ for which inequality \eqref{GrindEQ__90_} holds for all $\varphi \in \bar{S}$, $u\in H_{0}^{1} (0,1)$, $\left(\xi ,w\right)\in {\mathbb R}^{2} $. Combining \eqref{GrindEQ__90_} with \eqref{GrindEQ__116_} and using the fact that $n^{2} \pi ^{2} \left(\mho \bar{\varphi }_{0} ,g_{n} \right)^{2} =\left(\bar{\varphi }_{0} ,\phi _{n} \right)^{2} $ we obtain for $t\ge 0$: 
\begin{align} \label{GrindEQ__117_} 
&\frac{1}{2} \sum _{n=1}^{\infty }n^{2} \pi ^{2} \left(u[t],g_{n} \right)^{2}  +\sum _{n=1}^{\infty }n^{2} \pi ^{2} \left(1+\sigma n^{2} \pi ^{2} \right)\left(\varphi [t],\phi _{n} \right)^{2}  
\nonumber \\ 
\le& C\exp \left(-\omega t\right)\left(\sum _{n=1}^{\infty }\left(\bar{\varphi }_{0} ,\phi _{n} \right)^{2}  +\sum _{n=1}^{\infty }\left(1+\sigma n^{2} \pi ^{2} \right)n^{2} \pi ^{2} \left(\varphi _{0} ,\phi _{n} \right)^{2}  \right) \nonumber \\ &+\left(\frac{2}{\mu } +\frac{\mu +\kappa }{2\sigma } \right)\frac{8c_{5} c_{7} }{e\omega } \left(\sum _{n=1}^{\infty }\frac{1}{n^{2} \pi ^{2} }  \right)\exp \left(-\omega t\right) \nonumber \\ 
&\times \left(\xi _{0}^{2} +w_{0}^{2} +\left\| \mho \bar{\varphi }_{0} \right\| _{2}^{2} +\left\| \varphi _{0} \right\| _{H^{1} (0,1)}^{2} \right)  
\end{align}
Using \eqref{GrindEQ__104_}, the fact that $\pi ^{2} \left\| \mho \bar{\varphi }_{0} \right\| _{2}^{2} \le \left\| \bar{\varphi }_{0} \right\| _{2}^{2} $ (Wirtinger's inequality; recall definition \eqref{GrindEQ__8_}) and Parseval's identity we get from \eqref{GrindEQ__117_} for all $t\ge 0$:
\begin{align} \label{GrindEQ__118_} 
&\frac{1}{2} \left\| \varphi _{t} [t]\right\| _{2}^{2} +\left\| \varphi _{x} [t]\right\| _{2}^{2} +\sigma \left\| \varphi _{xx} [t]\right\| _{2}^{2}  \nonumber \\ 
\le& C\exp \left(-\omega t\right)\left(\left\| \bar{\varphi }_{0} \right\| _{2}^{2} +\left\| \varphi '_{0} \right\| _{2}^{2} +\sigma \left\| \varphi ''_{0} \right\| _{2}^{2} \right) \nonumber \\ 
&+\left(\frac{2}{\mu } +\frac{\mu +\kappa }{2\sigma } \right)\frac{8c_{5} c_{7} }{e\omega } \left(\sum _{n=1}^{\infty }\frac{1}{n^{2} \pi ^{2} }  \right)\exp \left(-\omega t\right) \nonumber \\
&\times \left(\xi _{0}^{2} +w_{0}^{2} +\frac{1}{\pi ^{2} } \left\| \bar{\varphi }_{0} \right\| _{2}^{2} +\left\| \varphi _{0} \right\| _{H^{1} (0,1)}^{2} \right)  
\end{align} 
Estimate \eqref{GrindEQ__103_} for certain appropriate constant $G>0$ is a consequence of \eqref{GrindEQ__61_}, \eqref{GrindEQ__118_} and the fact that $\pi ^{2} \left\| \mho \bar{\varphi }_{0} \right\| _{2}^{2} \le \left\| \bar{\varphi }_{0} \right\| _{2}^{2} $. The proof is complete.   
\qed
\end{proof}
\section{Conclusions}
The present paper intends to be a benchmark paper that will be followed by many other works on inverse optimality of PDE systems with boundary control. In this paper we managed to give a step-by-step procedure that allows the construction of families of optimal feedback laws without any restrictions coming from compatibility (or other avoidable regularity) conditions. It is clear that the present work provided results for a specific tank-liquid system. However, the procedure can be extended to various other PDE systems with boundary control. Future work will address such problems. 



\end{document}